\documentclass[11pt]{article}
\usepackage{latexsym}
\usepackage[all]{xy}
\usepackage{amsmath}
\usepackage{amssymb}
\usepackage{amsfonts}

\newtheorem{thm}{Theorem}[subsection]

\newtheorem{lem}[thm]{Lemma}
\newtheorem{df}[thm]{Definition}
\newtheorem{cor}[thm]{Corollary}

\newtheorem{rmk}[thm]{Remark}
\newtheorem{ex}[thm]{Example}
\newtheorem{Q}[thm]{Question}

\begin{document}

\title{\textbf{Brave New Algebraic Geometry\\ and global derived moduli spaces\\ of ring spectra}}
\bigskip
\bigskip

\author{\bigskip\\
Bertrand To\"en \\
\small{Laboratoire Emile Picard}\\
\small{UMR CNRS 5580} \\
\small{Universit\'{e} Paul Sabatier, Toulouse}\\
\small{France}\\
\bigskip \and
\bigskip \\
Gabriele Vezzosi \\
\small{Dipartimento di Matematica Applicata}\\
\small{``G. Sansone''}\\
\small{Universit\`a di Firenze}\\
\small{Italy}\\
\bigskip}

\maketitle

\begin{abstract}
We develop homotopical algebraic geometry (\cite{hagI,hagII}) in
the special context where the base symmetric monoidal model
category is that of spectra $\mathbf{S}$, i.e. what might be
called, after Waldhausen, \textit{brave new algebraic geometry}.
We discuss various model topologies on the model category of
commutative algebras in $\mathbf{S}$, and their associated
theories of geometric $\mathbf{S}$-stacks (a geometric
$\mathbf{S}$-stack being an analog of Artin notion of algebraic
stack in Algebraic Geometry). Two examples of geometric
$\mathbf{S}$-stacks are given: a global moduli space of
associative ring spectrum structures, and the stack of elliptic
curves endowed with the sheaf of topological modular forms.
\end{abstract}


\medskip

\textsf{MSC-class:} $55$P$43$; $14$A$20$; $18$G$55$; $55$U$40$;
$18$F$10$.


\tableofcontents


\setcounter{section}{0}

\bigskip
\bigskip
\bigskip

\begin{section}{Introduction}
Homotopical Algebraic Geometry is a kind of algebraic
geometry where the affine objects are given by commutative ring-like
objects in some homotopy theory (technically speaking, in a
symmetric monoidal model category); these affine objects are
then glued together according to an appropriate homotopical
modification of a Grothendieck topology (a \textit{model topology}, see \cite[4.3]{hagI}). More generally,
we allow ourselves to consider more flexible objects like stacks,
in order to deal with appropriate moduli problems.
This theory is developed in full generality in \cite{hagI, hagII} (see also \cite{hagdag}).
Our motivations for such a theory came from a variety of sources: first of
all, on the algebro-geometric side, we wanted to produce a sufficiently
functorial language in which the so called Derived Moduli Spaces foreseen
by Deligne, Drinfel'd and Kontsevich could really be constructed; secondly, on the topological side,
we thought that maybe the many recent results in Brave New Algebra,
i.e. in (commutative) algebra over structured ring spectra
(in any one of their brave new symmetric monoidal model categories,
see e.g. \cite{hss, ekmm}), could be pushed to a kind of Brave New
Algebraic Geometry in which one could take advantage of the possibility
of gluing these brave new rings together into an actual geometric object,
much in the same way as commutative algebra is helped (and generalized) by
the existence of algebraic geometry. Thirdly, on the motivic side, following
a suggestion of Y. Manin, we wished to have a sufficiently general theory in
order to study algebraic geometry \textit{over} the recent model categories of
motives for smooth schemes over a field (\cite{hu, ja3, sp}).\\

The purpose of this paper is to present the first steps in
the second type of applications mentioned above, i.e. a
specialization of the general framework of homotopical algebraic geometry
to the context of stable homotopy theory.
Our category $\mathbf{S-Aff}$ of brave new affine objects will therefore be defined as the
the opposite model category of the category of commutative rings in the category
$\mathbf{S}$ of symmetric spectra (\cite{hss}).

We first define and study various model topologies defined on $\mathbf{S-Aff}$.
They are all extensions, to different extents, of the usual Grothendieck topologies
defined on the category of (affine) schemes, like the Zariski and \'etale ones.

With any of these model topologies $\tau$ at our disposal, we
define and give the basic properties of the corresponding model
category of $\mathbf{S}$-stacks, understood in the broadest sense
as not necessarily truncated presheaves of simplicial sets on
$\mathbf{S-Aff}$ satisfying a homotopical descent (i.e.
sheaf-like) condition with respect to $\tau$-(hyper)covers. A
model topology on $\mathbf{S-Aff}$ is said to be
\textit{subcanonical} if the representable simplicial presheaves,
i.e. those of the form $\mathrm{Map}(A,-)$, for some commutative
ring $A$ in $\mathbf{S}$, $\mathrm{Map}$ being the mapping space
in $\mathbf{S-Aff}^{op}$, are $\mathbf{S}$-stacks.

As in algebraic geometry one finds it useful to study those stacks
defined by smooth groupoids (these are called Artin algebraic
stacks), we also define a brave new analog of these and call them
geometric $\mathbf{S}$-stacks, to emphasize that such
$\mathbf{S}$-stacks host a rich geometry very close to the
geometric intuition learned in algebraic geometry. In particular,
given a geometric $\mathbf{S}$-stack $F$, it makes sense to speak
about quasi-coherent and perfect modules over $F$, about the
$K$-theory of $F$, etc.; various properties of morphisms (e.g.
smooth, \'etale, proper, etc.) between geometric
$\mathbf{S}$-stacks can likewise be defined.

Stacks were introduced in algebraic geometry mainly to study moduli problems
of various sorts; they provide actual geometric objects (rather than sets of
isomorphisms classes or coarse moduli schemes) which store all the fine details
of the classification problem and on which a geometry very similar to that of
algebraic varieties or schemes can be developed, the two aspects having a fruitful
interplay. In a similar vein, in our brave new context, we give one example of a
moduli problem arising in algebraic topology (the classification of $A_{\infty}$-ring
spectrum structures on a given spectrum $M$) that can be studied geometrically
through the geometric $\mathbf{S}$-stack $\mathbb{R}\underline{Ass}_{M}$ it represents.
We wish to emphasize that instead of a discrete homotopy type (like the ones studied,
for different moduli problems, in \cite{re,bdg,gh}), we get a full geometric
object on which a lot of interesting geometry can be performed. The geometricity
of the $S$-stack $\mathbb{R}\underline{Ass}_{M}$, with respect to any fixed subcanonical model
topology, is actually the main theorem of this paper
(see Theorem \ref{t2}).

We also wish to remark that the approach presented in this paper can be extended to other, more interested and involved,
moduli problems algebraic topologists are interested in, and perhaps this richer
geometry could be of some help in answering, or at least in formulating in a clearer
way, some of the deep questions raised by the recent progress in stable homotopy theory (see \cite{goe}).
In this direction, we will explain in \S \ref{tmf} how
topological modular forms give rise to a natural
geometric $\mathbf{S}$-stack which is an extension in the brave new direction
of the moduli stack of elliptic curves (see Theorem \ref{t3}). This fact seems to us a very
important remark (probably much more interesting
than our Theorem \ref{t2}), and we think it could be the starting point of a very
interesting research program.

We also present a brave new analog of the stack of vector bundles
on a scheme, called the $\mathbf{S}$-(pre)stack $\mathsf{Perf}$ of perfect modules (Section 3.2), and
we expect it to be a key tool in brave new algebraic geometry.
The prestack $\mathsf{Perf}$ is a stack if and only if the model topology we are working with
is subcanonical (Thm. \ref{t1} whose proof is postponed to \cite{hagII}). This is another instance
of the relevance of the \textit{descent problem}, i.e. the question whether a given model topology
is subcanonical or not (see Section 3.1). Though we prove that some of the model topologies
we introduce (namely the standard and the semi-standard ones, Section 2.3) are subcanonical, at present we are
not able to settle (nor in the positive nor in the negative) the descent problem for the three
most promising model topologies we define, namely the Zariski, \'etale and thh-\'etale ones.
Though this is at the moment quite unsatisfactory, we believe that the descent problem
for these topologies is a very interesting question in itself even leaving outside its crucial role
in brave new algebraic geometry.

For the Zariski model topology, we have a partial positive result in this direction.
 By definition, for a model topology $\tau$ the property of being subcanonical depends on the notion of stacks
we consider; if instead of defining a stack as a prestack (i.e. a simplicial presheaf) satisfying
homotopical descent with respect to \textit{all} homotopy $\tau$-hypercovers, we simply require
descent with respect to all $\check{\mathrm{C}}$ech $\tau$-hypercovers (i.e. those arising
as nerves of $\tau$-covers), we obtain a notion of \textit{$\check{C}$ech stacks}, recently
considered by J. Lurie (\cite{lu}) and Dugger-Hollander-Isaksen (\cite{dhi}). We prove (Corollary \ref{zariscechsubcan})
that the Zariski model topology is in fact subcanonical with respect to the notion of $\check{\mathrm{C}}$ech stacks.
Moreover, by replacing in Theorem \ref{t2} the word ``stack'' with the weaker
``$\check{\mathrm{C}}$ech stack'', the statement remains true for any model topology which is
subcanonical with respect to the notion of $\check{\mathrm{C}}$ech stacks.

It is therefore natural to ask why we did not choose to formulate everything only in terms of
$\check{\mathrm{C}}$ech stacks. We believe that at this early stage of development of homotopical
algebraic geometry and, in particular, of brave new algebraic geometry, it is not advisable to
make choices that could prevent some applications
or obscure some of the properties of the objects involved, while it is more useful to keep in mind
various options, some of which can be more useful in one context than in others. For example,
it is clear that knowing that a given, geometrically meaningful, simplicial presheaf is a stack and not
only a $\check{\mathrm{C}}$ech stack adds a lot more informations, in fact exactly the descent property
with respect to unbounded hypercovers (\cite[Thm. A.6]{dhi}). Moreover, $\check{\mathrm{C}}$ech stacks fail
in general to satisfy an analog of Whitehead theorem: a pointed $\check{\mathrm{C}}$ech stack
may have vanishing $\pi_{i}$ sheaves for any $i\geq 0$ without being necessarily contractible. This
last fact is a very inconvenient property of $\check{\mathrm{C}}$ech stacks, that
makes Postnikov decompositions and spectral sequences arguments uncertain.
On the other hand, the $\check{\mathrm{C}}$ech descent condition is usually much easier to establish
than the full descent condition, and as we have already remarked, some natural model topologies
are easily seen to be subcanonical with respect to the notion of $\check{\mathrm{C}}$ech stacks while
it might be tricky to show that they are actually subcanonical. Finally, we would like to
mention that in our experience we have never met
serious troubles by using one or the other of the two notions, and
in many interesting contexts it does not really matter which notion one uses,
as the rather subtle differences actually tend not to appear in practice. \\

\textbf{Acknowledgments}. We wish to thank the organizers of the INI Workshop on Elliptic
Cohomology and Higher Chromatic Phenomena (Cambridge UK, December 2003) for the invitation
to speak about our work
and Bill Dwyer for his encouraging comments. We would also like
to thank Michael Mandell, Peter May, Haynes Miller,
John Rognes, Stefan Schwede and Neil Strickland for
helpful discussions and suggestions.

\end{section}

\bigskip

\textbf{Notations.}
To fix ideas, we will work in the category
$\mathbf{S}:=\mathbf{Sp}^{\Sigma}$ of symmetric spectra
(see \cite{hss}), but all the constructions of this paper will also work,
possibly with minor variations (see \cite{schw}), for other equivalent theories (e.g for the category of
$S$-modules of \cite{ekmm}). We will consider $\mathbf{S}$ as a
symmetric monoidal simplicial model category (for the smash product $-\wedge -$)
with the Shipley-Smith positive $S$-model
structure (see \cite[Prop. 3.1]{shi}).

We define $\mathbf{S-Alg}$ as the category of (associative and
unital) commutative monoids objects in $\mathbf{S}$, endowed with
the $S$-model structure of \cite[Thm. 3.2]{shi}; we will simply
call them \emph{commutative $S$-algebras} instead of the more correct but longer, commutative
symmetric ring spectra. For any commutative $S$-algebra $A$,
we will denote by $A-\mathbf{Alg}$ the under-category
$A/\mathbf{S-Alg}$, whose objects will be called \emph{commutative
$A$-algebras}. Finally, if $A$ is a commutative $S$-algebra,
$A-\mathbf{Mod}$ will be the category of $A$-modules with the
$A$-model structure (\cite[Prop. 3.1]{shi}). This model category
is also a symmetric monoidal model category for the smash product $-\wedge_{A} -$
over $A$.

For a morphism of commutative $S$-algebras, $f : A \longrightarrow B$
one has a Quillen adjunction
$$f^{*} : A-\mathbf{Mod} \longrightarrow B-\mathbf{Mod}
\qquad A-\mathbf{Mod} \longleftarrow B-\mathbf{Mod} : f_{*},$$
where $f^{*}(-):=-\wedge_{A}B$ is the base change functor.
We will denote by
$$\mathbb{L}f^{*} : \mathbf{Ho}(A-\mathbf{Mod}) \longrightarrow \mathbf{Ho}(B-\mathbf{Mod})
\qquad \mathbf{Ho}(A-\mathbf{Mod}) \longleftarrow \mathbf{Ho}(B-\mathbf{Mod}) : \mathbb{R}f_{*}$$
the induced derived adjunction on the homotopy categories.

Our references for model category theory are
\cite{hi,ho}. For a model category $M$ with equivalences $W$, the set of morphisms
in the homotopy category $\mathbf{Ho}(M):=W^{-1}M$ will be
denoted by $[-,-]_{M}$, or simply by
$[-,-]$ if the context is clear. The (homotopy) mapping spaces in $M$ will be denotedby $\mathrm{Map}_{M}(-,-)$. When $M$ is a simplicial model category,
the simplicial Hom's (resp. derived simplicial Hom's) will be denoted by
$\underline{Hom}_{M}$ (resp. $\mathbb{R}\underline{Hom}_{M}$), or simply by
$\underline{Hom}$ (resp. $\mathbb{R}\underline{Hom}$) if the context is clear.
Recall that in this case one can compute  $\mathrm{Map}_{M}(-,-)$ as $\mathbb{R}\underline{Hom}_{M}(-,-)$.

Finally, for a model category $M$ and an object $x\in M$ we will often use the
coma model categories $x/M$ and $M/x$. When the model category $M$ is not
left proper (resp. is not right proper) we will always assume that
$x$ has been replaced by a cofibrant (resp. fibrant)
model before considering $x/M$ (resp. $M/x$). More generally,
we will not always mention fibrant and cofibrant replacements and suppose implicitly that
all our objects are fibrant and/or cofibrant when required. \\

\bigskip

Since we wish to concentrate on applications
to stable homotopy theory, some general constructions and
details about homotopical algebraic geometry will be
omitted by referring to \cite{hagI}.
For a few of the results presented we will only give here sketchy proofs; full proofs will appear in \cite{hagII}.\\

\bigskip

\begin{section}{Brave new sites}

In this section we present two model topologies on the (opposite) category
of commutative $S$-algebras. They are \emph{brave new analogs} of the
Zariski and \'etale topologies defined on the category of usual commutative
rings and will allow us to define the \emph{brave new Zariski} and \emph{\'etale sites}. \\

We denote by $\mathbf{S-Aff}$ the opposite model category of
$\mathbf{S-Alg}$. \\

If $M$ is a model category we say that an object $x$ in $M$ is \textit{finitely presented} if, for any filtered direct system of objects $\left\{z_{i}\right\}_{i\in J}$ in $M$, the natural map $$\mathrm{colim}_{i}\; \mathrm{Map}_{M}(x,z_{i})\longrightarrow \mathrm{Map}_{M}(x,\mathrm{colim}_{i}\;z_{i})$$ is an isomorphism in the homotopy category of simplicial sets. \\

\begin{df}\label{d0}
A morphism $A\rightarrow B$ of commutative $S$-algebras is \emph{finitely presented} if it is a finitely presented object in the model under-category $A/(\mathbf{S-Alg})= A-\mathbf{Alg}$; in this case, we also say that $B$ is
a finitely presented commutative $A$-algebra. An $A$-module $E$ is \emph{finitely presented} or \emph{perfect} if it is a finitely presented object in the model category $A-\mathbf{Mod}$. \\
\end{df}

Perfect $A$-modules can also be characterized as retracts of finite cell $A$-modules (see \cite[Thm. III-7.9]{ekmm}); in particular, there are plenty of them. If $A$ is a commutative $S$-algebra, then the free
commutative $A$-algebra on a finite number of generators (or, more generally, on any perfect $A$-module) is a finitely presented $A$-algebra.
The reader will find other examples of finitely presented morphisms
of commutative $S$-algebras in Lemma \ref{l2}.

\begin{subsection}{The brave new Zariski topology}

\begin{df}\label{d1}
\begin{itemize}
\item A morphism $f : A \longrightarrow B$ in $\mathbf{S-Alg}$ is called
a \emph{formal Zariski open immersion} if the induced functor
$\mathbb{R}f_{*} : \mathbf{Ho}(B-\mathbf{Mod})
\longrightarrow \mathbf{Ho}(A-\mathbf{Mod})$ is fully faithful.
\item A morphism $f : A\longrightarrow B$ is a \emph{Zariski open immersion} if
$\mathbf{S-Alg}$ is
it is a formal Zariski open immersion and
of finite presentation (as a morphism of commutative
$S$-algebras).
\item A family $\{f_{i} : A \longrightarrow A_{i}\}_{i\in I}$ of morphisms
in $\mathbf{S-Alg}$ is called a
\emph{(formal) Zariski open covering} if it satisfies the following two conditions.
\begin{itemize}

\item Each morphism $A \longrightarrow A_{i}$ is a (formal) Zariski open immersion.

\item There exist a finite subset $J \subset I$ such that the
family of inverse image functors
$$\{\mathbb{L}f_{j}^{*} : \mathbf{Ho}(A-\mathbf{Mod}) \longrightarrow \mathbf{Ho}(A_{j}-\mathbf{Mod})\}_{j\in J}$$
is conservative (i.e. a morphism in $\mathbf{Ho}(A-\mathbf{Mod})$ is an isomorphism if and only if
its images by all the $\mathbb{L}f_{j}^{*}$'s are isomorphisms).
\end{itemize}
\end{itemize}
\end{df}

\begin{ex}\emph{If $A\in \mathbf{S-Alg}$ and $E$ is an $A$-module such
that the associated Bousfield localization $L_{E}$ is smashing
(i.e. the natural transformation $L_{E}(-)\rightarrow L_{E}A\wedge^{\mathbb{L}}_{A}(-)$
is an isomorphism), then $A\rightarrow L_{E}A$ (which is a morphism
of commutative $S$-algebras by e.g. \cite[\S VIII.2]{ekmm}) is a formal
Zariski open immersion. This follows immdiately from the fact that
$\mathrm{\mathbf{Ho}}(L_{E}A-Mod)$ is equivalent to the subcategory of $\mathrm{\mathbf{Ho}}(A-Mod)$
consisting of $L_{E}$-local objects, by \cite{wo}.}
\end{ex}

It is easy to check that (formal) Zariski open covering families define
a \textit{model topology} in the sense of \cite[\S 4.3]{hagI} on the model category $\mathbf{S-Aff}$. For the reader's
convenience we recall what this means in the following lemma.

\begin{lem}\label{l0}
\begin{itemize}

\item If $A \longrightarrow B$ is an equivalence of commutative $S$-algebras then the
one element family $\{A \longrightarrow B\}$ is a (formal) Zariski open covering.

\item Let $\{A \longrightarrow A_{i}\}_{i\in I}$ be a (formal) Zariski open covering of
$S$-algebras and $A \longrightarrow B$ a morphism. Then, the family of homotopy push-outs
$\{B \longrightarrow B\wedge_{A}^{\mathbb{L}}A_{i}\}_{i\in I}$ is also a (formal) Zariski open
covering.

\item Let $\{A \longrightarrow A_{i}\}_{i\in I}$ be a (formal) Zariski open covering of
$S$-algebras, and for any $i \in I$ let
$\{A_{i} \longrightarrow A_{ij}\}_{j\in J_{i}}$ be a (formal) Zariski open covering of
$S$-algebras. Then, the total family
$\{A \longrightarrow A_{ij}\}_{i\in I,j\in J_{i}}$ is again a (formal) Zariski open covering.
\end{itemize}
\end{lem}

\textit{Proof:} Left as an exercise. \hfill $\Box$ \\

By definition,
Lemma \ref{l0} shows that (formal) Zariski open coverings define
a model topology on the model category $\mathbf{S-Aff}$ and so, as proved in \cite[Prop. 4.3.5]{hagI}, induce a Grothendieck
topology on the homotopy category $\mathbf{Ho}(\mathbf{S-Alg})$.
This model topology is called the \emph{brave new (formal) Zariski topology}, and
endows $\mathbf{S-Aff}$ with the structure of a model site in the sense
of \cite[\S 4]{hagI}. This model site, denoted by $(\mathbf{S-Aff},\mathrm{Zar})$ for the brave new Zariski topology, and
$(\mathbf{S-Aff},\mathrm{fZar})$ for the brave new formal Zariski topology. They will be called
the \emph{brave new Zariski site} and the \emph{brave new formal Zariski site}. \\

Let $\mathbf{Alg}$ be the category of (associative and unital) commutative
rings.
Let us recall the existence of the Eilenberg-Mac Lane functor
$$H : \mathbf{Alg} \longrightarrow \mathbf{S-Alg},$$
sending a commutative ring $R$ to the commutative $S$-algebra $HR$
such that $\pi_{0}(HR)=R$ and $\pi_{i}(HR)=0$ for any $i\neq 0$.
This functor is homotopically fully faithful and
the following lemma shows that our brave new Zariski topology does
generalize the usual Zariski topology.

\begin{lem}\label{l1}
\begin{enumerate}
\item[\emph{(1)}]
Let $R \longrightarrow R'$ be a morphism of commutative rings.
The induced morphism $HR \longrightarrow HR'$ is a Zariski open immersion
of commutative $S$-algebras (in the sense of Definition \ref{d1}) if and only if
the morphism $Spec\, R' \longrightarrow Spec\, R$ is an open immersion of schemes.

\item[\emph{(2)}]
A family of morphisms of
commutative rings, $\{R \longrightarrow R'_{i}\}_{i\in I}$, induces
a Zariski covering family of commutative $S$-algebras $\{HR \longrightarrow HR_{i}'\}_{i\in I}$
(in the sense of Definition \ref{d1})
if and only if the family $\{Spec\, R_{i} \longrightarrow Spec\, R\}_{i\in I}$
is a Zariski open covering of schemes.
\end{enumerate}
\end{lem}

\textit{Proof:} Let us start with the general situation of a morphism $f : A \longrightarrow B$
of commutative $S$-algebras such that
the induced functor $\mathbb{R}f_{*} : \mathbf{Ho}(B-\mathbf{Mod}) \longrightarrow \mathbf{Ho}(A-\mathbf{Mod})$
is fully faithful. Let $L=\mathbb{R}f_{*}\circ \mathbb{L}f^{*}$, which comes
with a natural transformation $Id \longrightarrow L$.
Then, the essential image of $\mathbb{R}f_{*}$ consist of objects $M$
in $\mathbf{Ho}(A-\mathbf{Mod})$ such that the localization morphism $M \longrightarrow LM$
is an isomorphism.

The Quillen adjunction $(f^{*},f_{*})$
extends to a Quillen adjunction on the category of commutative algebras
$$f^{*} : A-\mathbf{Alg} \longrightarrow B-\mathbf{Alg} \qquad
A-\mathbf{Alg} \longleftarrow B-\mathbf{Alg} : f_{*},$$ also with the property that $\mathbb{R}f_{*} : \mathbf{Ho}(B-\mathbf{Alg})
\longrightarrow \mathbf{Ho}(A-\mathbf{Alg})$ is fully faithful.
Furthermore, the essential image of this last functor consist of
all objects $C \in \mathbf{Ho}(A-\mathbf{Alg})$ such that the underlying
$A$-module of $C$ satisfies $C\simeq LC$ (i.e. the underlying $A$-module
of $C$ lives in the image of $\mathbf{Ho}(B-\mathbf{Mod})$).

>From these observations, we deduce that for any commutative $A$-algebra
$C$, the mapping space $\mathbb{R}\underline{Hom}_{A-\mathbf{Alg}}(B,C)$
is either empty or contractible; it is non-empty if and only if
the underlying $A$-module of $C$ belongs to the essential image of $\mathbb{R}f_{*}$. \\

To prove $(1)$,
let us first suppose that $f : Spec\, R'  \longrightarrow Spec\, R$
is an open immersion of schemes. The induced functor
on the derived categories $f_{*} : D(R') \longrightarrow D(R)$ is then fully faithful.
As there are natural equivalences (\cite[IV Thm. 2.4]{ekmm})
$$\mathbf{Ho}(HR-\mathbf{Mod})\simeq D(R) \qquad \mathbf{Ho}(HR'-\mathbf{Mod})\simeq D(R')$$
this implies that the functor
$\mathbb{R}f_{*} : \mathbf{Ho}(HR'-\mathbf{Mod}) \longrightarrow \mathbf{Ho}(HR-\mathbf{Mod})$ is also
fully faithful. It only remains to show that $HR \longrightarrow HR'$ is
finitely presented in the sense of Definition \ref{d0}.

We will first assume that $R'=R_{f}$ for some element $f\in R$.
The essential image of $\mathbb{R}f_{*} : \mathbf{Ho}(HR'-\mathbf{Mod}) \longrightarrow \mathbf{Ho}(HR-\mathbf{Mod})$
then consists of all objects $E \in \mathbf{Ho}(HR-\mathbf{Mod})\simeq D(R)$
such that $f$ acts by isomorphisms on the cohomology $R$-module $H^{*}(E)$.
By what we have seen at the beginning of the proof,
this implies that for any commutative $HR$-algebra $C$
the mapping space $\mathbb{R}\underline{Hom}_{HR-\mathbf{Alg}}(HR',C)$ is
contractible if $f$ becomes invertible in $\pi_{0}(C)$, and empty otherwise.
>From this one easily deduces that
$\mathbb{R}\underline{Hom}_{HR-\mathbf{Alg}}(HR',-)$ commutes with filtered colimits, or in other words that
$HR'$ is a finitely presented $HR$-algebra in the sense of Definition \ref{d0}.

In the general case, one can write $Spec\, R'$ as a finite union
of schemes of the form $Spec\, R_{f}$ for some elements $f\in R$. A bit of descent theory (see
\S 3.1) then allows us to reduce to the case where $R'=R_{f}$ and conclude. \\

Let us now assume that $HR \longrightarrow HR'$ is
a Zariski open immersion of commutative $S$-algebras. By adjunction (between
$H$ and $\pi_{0}$ restricted on connective $S$-algebras) one sees easily
that $R \longrightarrow R'$ is a finitely presented morphism of
commutative rings.

The induced
functor on (unbounded) derived categories
$$f_{*} : D(R')\simeq \mathbf{Ho}(HR'-\mathbf{Mod}) \longrightarrow D(R)\simeq \mathbf{Ho}(HR-\mathbf{Mod})$$
is fully faithful. Through the Dold-Kan correspondence, this implies that
the Quillen adjunction on the model category
of simplicial modules (see \cite{gj})
$$f^{*} : sR-\mathbf{Mod} \longrightarrow sR'-\mathbf{Mod} \qquad sR-\mathbf{Mod} \longleftarrow
sR'-\mathbf{Mod} : f_{*}$$
is such that $\mathbb{L}f^{*}\circ f_{*}\simeq Id$.
Let $sR-\mathbf{Alg}$ and $sR'-\mathbf{Alg}$ be the categories of simplicial
commutative $R$-algebras and simplicial commutative $R'$-algebras, endowed
with their natural model structures (equivalences are and fibration are
detected in the category of simplicial modules). Then, the Quillen adjunction
$$f_{*} : sR'-\mathbf{Alg} \longrightarrow sR-\mathbf{Alg} \qquad sR'-\mathbf{Alg} \longleftarrow
sR-\mathbf{Alg} : f^{*}$$
also satisfies $\mathbb{L}f^{*}\circ f_{*}\simeq Id$, as this is true on the level on
simplicial modules. In particular, for any simplicial $R'$-module $M$,
the space of derived derivations
$$\mathbb{L}Der_{R}(R',M):=\mathbb{R}\underline{Hom}_{sR-\mathbf{\mathbf{Alg}}/R'}(R',R'\oplus M)\simeq
\mathbb{R}\underline{Hom}_{sR'-\mathbf{Alg}/R'}(R',R'\oplus M)\simeq *$$
is acyclic (here $R'\oplus M$ is the simplicial $R'$-algebra which is the trivial extension
of $R'$ by $M$). As a consequence one sees that
 Quillen's cotangent complex $\mathbb{L}_{R'/R}$ is acyclic, which implies that the
morphism $R \longrightarrow R'$ is an \'etale morphism of rings.

Finally, using the fact that the functor on the category of modules $R'-\mathbf{Mod} \longrightarrow
R-\mathbf{Mod}$
is fully faithful, one sees that $Spec\, R' \longrightarrow Spec\, R$ is
a monomorphism of schemes. Therefore, the morphism of schemes
$Spec\, R' \longrightarrow Spec\, R$ is an \'etale monomorphism, and so
is an open immersion by \cite[Thm. 17.9.1]{ega4}. \\

Finally, point $(2)$ is clear if one knows $(1)$ and that
$\mathbf{Ho}(HR-\mathbf{Mod})\simeq D(R)$. \hfill $\Box$ \\

\begin{rmk}\label{loc}
\emph{The argument at the beginning of the proof of Lemma \ref{l1} shows that if $f:A\rightarrow B$ is a Zariski open immersion, the functor $L(f):=\mathbb{R}f_{*}\mathbb{L}f^{*}$ is a \textit{localization functor} on the homotopy category of $A$-modules in the sense of \cite[Def. 3.1.1]{hps}. And it is also clear by definition that $L(f)$ is also \textit{smashing} (\cite[Def. 3.3.2]{hps}). Let us call a localization functor $L$ on $\mathbf{Ho}(A-Mod)$ a \textit{formal Zariski localization functor over} $A$ if $L\simeq L(f)$ for some formal Zariski open immersion $f$. Let us also say that a localization functor $L$ on $\mathrm{\mathbf{Ho}}(A-Mod)$ is a \textit{smashing algebra Bousfield localization over $A$} if $L\simeq L_{B}$ for some $A$-algebra $B$ such that $L_{B}$ is smashing (over $A$). Then it is easy to verify that in the set of equivalence classes of localization functors on $\mathrm{\mathbf{Ho}}(A-Mod)$, the subset consisting of formal Zariski localization functors over $A$ coincides with the subset consisting of smashing algebra Bousfield localizations over $A$. In fact, if $f:A\rightarrow B$ is a Zariski open immersion, $L_{B}$ denotes the Bousfield localization with respect to the $A$-module $B$, and $\ell_{B/A}:A\rightarrow L_{B}A$ the corresponding morphism of commutative $A$-algebras, we have $L(f)\simeq L_{B}\simeq L(\ell_{B/A})$ because all three localizations have the same category of acyclics. Conversely, if $L_{C}$ is a smashing algebra Bousfield localization over $A$, and $\ell_{C/A}:A\rightarrow L_{C}A$ is the corresponding morphism of commutative $A$-algebras, one has $L_{B}\simeq L(\ell_{C/A})$.  }
\end{rmk}

Let $\mathbf{Aff}$ be the opposite category of commutative rings, and
$(\mathbf{Aff},\mathrm{Zar})$ the big Zariski site. The site
$(\mathbf{Aff},\mathrm{Zar})$ can also be considered as a model site (for the trivial
model structure on $\mathbf{Aff}$). Lemma \ref{l1} implies in particular that
the functor $H : \mathbf{Aff} \longrightarrow \mathbf{S-Aff}$
induces a continuous morphism of model sites (\cite[Def. 4.8.4]{hagI}). In this way, the
site $(\mathbf{Aff},\mathrm{Zar})$ becomes a \emph{sub-model site}
of $(\mathbf{S-Aff},\mathrm{Zar})$. \\

To finish with the Zariski topology we will now describe a general
procedure in order to construct interesting open Zariski immersions of commutative $S$-algebras
using the techniques of Bousfield localization for model categories. \\

Let $A$ be a commutative $S$-algebra and $M$ be a $A$-module. We will assume that
$M$ is a perfect $A$-module (in the sense of Definition \ref{d0}), or equivalently that it is a strongly dualizable object in the monoidal category
$\mathbf{Ho}(A-\mathbf{Mod})$. As already noticed, perfect $A$-modules
are exactly the retracts of finite cell $A$-modules, see \cite[Thm. III-7.9]{ekmm}).
Let $M[n]=S^{n}\otimes^{\mathbb{L}} M$ be the $n$-th suspension $A$-module of $M$, for
$n\in \mathbb{Z}$.

We denote by $D(M[n])$ the derived dual of $M[n]$, defined as
the derived internal Hom's of $A$-modules
$$D(M[n]):=\mathbb{R}\mathcal{HOM}_{A-\mathbf{Mod}}(M[n],A).$$
Consider now the (derived) free commutative $A$-algebra over $D(M[n])$,
$\mathbb{L}F_{A}(D(M[n]))$, characterized by the usual adjunction
$$[\mathbb{L}F_{A}(D(M[n])),-]_{A-\mathbf{Alg}}\simeq [D(M[n]),-]_{A-\mathbf{Mod}}.$$
The model category $A-\mathbf{Alg}$ is a combinatorial and
cellular model category, and therefore one can apply the
localization techniques (see e.g. \cite{hi,sm}) in order to invert
the natural augmentations $F_{A}(D(M[n])) \longrightarrow A$ for all $n\in \mathbb{Z}$.
One checks easily that, since
$M$ is strongly dualizable, the local objects for this
localization are the commutative $A$-algebras $B$ such that
$M\wedge_{A}^{\mathbb{L}}B\simeq 0$ in $\mathbf{Ho}(B-\mathbf{Mod})$.
The local model of $A$ for this localization will be denoted by
$A_{M}$. By definition, it is characterized by the following universal
property: for any commutative $A$-algebra $B$, the mapping space
$\mathbb{R}\underline{Hom}_{A-\mathbf{Alg}}(A_{M},B)$ is
contractible if $B\wedge^{\mathbb{L}}_{A}M\simeq 0$ and
empty otherwise. In other words, for any commutative $S$-algebra $B$ the
natural morphism
$$\mathbb{R}\underline{Hom}_{\mathbf{S-Alg}}(A_{M},B) \longrightarrow
\mathbb{R}\underline{Hom}_{\mathbf{S-Alg}}(A,B)$$ is equivalent to
an inclusion of connected components and its image consists of
morphisms $A \longrightarrow B$ in $\mathbf{Ho}(\mathbf{S-Alg})$ such that
$B\wedge^{\mathbb{L}}_{A}M\simeq 0$.

\begin{lem}\label{l2}
With the above notations, the morphism $A \longrightarrow A_{M}$ is a Zariski open immersion.
\end{lem}

\textit{Proof:} Let us start by showing that $A_{M}$ is a finitely presented
commutative $A$-algebra.

Let $\{B_{i}\}_{i\in I}$ be a filtered system of commutative $A$-algebras
and $B=\mathrm{colim}_{i}B_{i}$. We assume that $B\wedge^{\mathbb{L}}_{A}M\simeq 0$, and
we need to prove that there exists an $i\in I$ such that
$B_{i}\wedge_{A}^{\mathbb{L}}M\simeq 0$.

By assumption, the two points $Id$ and $0$ are the same in
$\pi_{0}(\mathbb{R}\underline{End}_{B-\mathbf{Mod}}(M\wedge_{A}^{\mathbb{L}}B))$. But, as
$M$ is a perfect $A$-module one has
$$\pi_{0}(\mathbb{R}\underline{End}_{B-\mathbf{Mod}}(M\wedge_{A}^{\mathbb{L}}B))\simeq
\mathrm{colim}_{i\in I}\pi_{0}(\mathbb{R}\underline{End}_{B_{i}-\mathbf{Mod}}(M\wedge_{A}^{\mathbb{L}}B_{i})).$$
This implies that there is some index $i\in I$ such that $Id$ and $0$ are
homotopic in $\mathbb{R}\underline{End}_{B_{i}-\mathbf{Mod}}(M\wedge_{A}^{\mathbb{L}}B_{i})$, and
therefore that $M\wedge_{A}^{\mathbb{L}}B_{i}$ is contractible. 

It remains to prove that the induced functor
$\mathbf{Ho}(A_{M}-\mathbf{Mod}) \longrightarrow \mathbf{Ho}(A-\mathbf{Mod})$ is fully faithful. For this, one uses that for any 
commutative $S$-algebra $B$, the morphism
$$\mathbb{R}\underline{Hom}_{\mathbf{S-Alg}}(A_{M},B) \longrightarrow
\mathbb{R}\underline{Hom}_{\mathbf{S-Alg}}(A,B)$$
is an inclusion of connected components. Therefore, the natural 
morphism 
$$\mathbb{R}\underline{Hom}_{\mathbf{S-Alg}}(A_{M},B) \longrightarrow
\mathbb{R}\underline{Hom}_{\mathbf{S-Alg}}(A_{M},B)\times_{
\mathbb{R}\underline{Hom}_{\mathbf{S-Alg}}(A,B)}^{h}
\mathbb{R}\underline{Hom}_{\mathbf{S-Alg}}(A_{M},B)$$
is an isomorphism in $Ho(SSet)$. This implies that 
the natural morphism
$$A_{M} \longrightarrow A_{M}\wedge^{\mathbb{L}}_{A}A_{M}$$
is an equivalence of commutative $S$-algebras. In particular, one
has for any $A_{M}$-module $M$
$$M\simeq M\wedge_{A_{M}}^{\mathbb{L}}A_{M}\simeq 
M\wedge_{A_{M}}^{\mathbb{L}}(A_{M}\wedge^{\mathbb{L}}_{A}A_{M})
\simeq M\wedge_{A}^{\mathbb{L}}A_{M},$$
showing that the base change functor
$$\mathbf{Ho}(A_{M}-\mathbf{Mod}) \longrightarrow \mathbf{Ho}(A-\mathbf{Mod})$$
is fully faithful.  \hfill $\Box$ \\

An important property of the localization $A\longrightarrow A_{M}$ is the following fact.

\begin{lem}\label{l3}
Let $A$ be a commutative $S$-algebra, and $M$ be a perfect $A$-module.
Then the essential image of the fully faithful functor
$$\mathbf{Ho}(A_{M}-\mathbf{Mod}) \longrightarrow \mathbf{Ho}(A-\mathbf{Mod})$$
consists of all $A$-modules $N$ such that
$M\wedge_{A}^{\mathbb{L}}N\simeq D(M)\wedge_{A}^{\mathbb{L}}N\simeq 0$.
\end{lem}

Note that since $M$ is perfect, then for any $A$-module $N$, $M\wedge_{A}^{\mathbb{L}}N\simeq 0$ iff  $D(M)\wedge_{A}^{\mathbb{L}}N\simeq 0$, so the two conditions in the lemma are actually one. Moreover, $A_{M}\simeq A_{D(M)}$ in $\mathbf{Ho}(A-\mathbf{Alg})$.\\

\textit{Proof:} As every $A_{M}$-module $N$ can be constructed by homotopy colimits of
free $A_{M}$-modules and $-\wedge_{A}^{\mathbb{L}}M$ commutes with
homotopy colimits, it is clear that $A_{M}\wedge^{\mathbb{L}}_{A}M\simeq 0$
implies $N\wedge_{A}^{\mathbb{L}}M\simeq 0$. Since
$A_{M}\wedge^{\mathbb{L}}_{A}D(M)\simeq D(A_{M}\wedge^{\mathbb{L}}_{A}M)\simeq 0$
(here the second derived dual is in the category of $A_{M}$-modules), the same argument shows that
$N\wedge_{A}^{\mathbb{L}}D(M)\simeq 0$.

Conversely, let $N$ be an $A$-module such that $N\wedge_{A}^{\mathbb{L}}M\simeq N\wedge^{\mathbb{L}}_{A}D(M)\simeq 0$.
By definition, the commutative $A$-algebra $A \longrightarrow A_{M}$ is obtained as a local model
of $A\rightarrow A$ when one inverts the set of morphisms
$\mathbb{L}F_{A}(D(M[n])) \longrightarrow A$, for any $n\in \mathbb{Z}$. It is well known (see e.g. \cite[\S 4]{hi}) that such a local model can be
obtained by a transfinite composition of homotopy push-outs of the form
$$\xymatrix{
A_{\alpha} \ar[r] & A_{\alpha+1} \\
\partial\Delta^{p}\otimes^{\mathbb{L}}\mathbb{L}F_{A}(D(M[n])) \ar[r] \ar[u] &
\Delta^{p}\otimes^{\mathbb{L}}\mathbb{L}F_{A}(D(M[n]))\ar[u]}$$ in the category of $A$-algebras.
>From this description, and the fact that $-\wedge_{A}^{\mathbb{L}}M$ commutes with
homotopy colimits, one sees that the adjunction morphism $N\longrightarrow N\wedge_{A}^{\mathbb{L}}A_{M}$
is an equivalence because by assumption on $N$, the natural morphism
$N\simeq N\wedge_{A}^{\mathbb{L}}A \longrightarrow N\wedge_{A}^{\mathbb{L}}\mathbb{L}F_{A}(D(M[n]))$ is an equivalence.
\hfill $\Box$ \\

Lemma \ref{l3} allows us to interpret geometrically $A_{M}$ as the \emph{open complement} of the
\emph{support} of the $A$-module $M$. Lemma \ref{l3} also has a converse whose proof is left as an exercise.

\begin{lem}\label{l3'}
Let $f : A \longrightarrow B$ be a morphism of commutative $S$-algebras and $M$ be a perfect
$A$-module. We suppose that the functor $\mathbb{R}f_* : \mathbf{Ho}(B-\mathbf{Mod}) \longrightarrow
\mathbf{Ho}(A-\mathbf{Mod})$ is fully faithful and that its essential image consists
of all $A$-modules $N$ such that $N\wedge_{A}^{\mathbb{L}}M\simeq 0$. Then,
the two commutative $A$-algebras $B$ and $A_{M}$ are equivalent (i.e. isomorphic
in $\mathbf{Ho}(A-\mathbf{Alg})$).
\end{lem}

\begin{rmk}\label{r1}
\emph{One should note carefully that even though if the Eilenberg-Mac Lane functor $H$ embeds $(\mathbf{Aff},\mathrm{Zar})$
in $(\mathbf{S-Aff},\mathrm{Zar})$ as model sites, there exist commutative rings
$R$ and Zariski open coverings $HR \longrightarrow B$ in
$\mathbf{S-Aff}$ such that $B$ is not of the form $HR'$ for some
commutative $R$-algebra $R'$. One example is given by taking $R$
to be $\mathbb{C}[X,Y]$, and considering the localized commutative
$HR$-algebra $(HR)_{M}$ (in the sense above), where $M$ is the perfect $R$-module
$R/(X,Y)\simeq \mathbb{C}$. If $(HR)_{M}$ were of the form $HR'$
for a Zariski open immersion $Spec\, R' \longrightarrow Spec\, R$,
then for any other commutative ring $R''$, the set of scheme morphisms
$Hom(Spec\, R'',Spec\, R')$ would be the subset of $Hom(Spec\,
R'',\mathbb{A}^{2})$ consisting of morphisms
factoring through $\mathbb{A}^{2}-\{0\}$. This would mean that
$Spec\, R'\simeq \mathbb{A}^{2}-\{0\}$, which is not possible as
$\mathbb{A}^{2}-\{0\}$ is a not an affine scheme. This example
is of course the same as the example given in \cite[\S 2.2]{to1} of a $0$-truncated affine stack
which is not an affine scheme. These kind of example shows that
there are many more affine objects in homotopical algebraic geometry than in
usual algebraic geometry.}
\end{rmk}

\begin{rmk}\label{nilpotence}
\emph{
\begin{enumerate}
  \item Note that Lemma \ref{l3} shows that the localization
process $\xymatrix{(A,M)\ar@{~>}[r] & A_{M}}$
is in some sense ``orthogonal'' to the usual
Bousfield localization process
$\xymatrix{(A,M)\ar@{~>}[r] & L_{M}A}$ in that the local objects for the
former are exactly the acyclic objects for
the latter. To state everything in terms of Bousfield localizations,
this says that $L_{A_{M}}$-local objects are
exactly $L_{M}$-acyclic objects (compare with Remark \ref{loc}).
Note that however, while the Bousfield localization is always defined for any $A$-module $M$,
the commutative $A$-algebra $A_{M}$ probably does not exist unless $M$ is perfect.
  \item Let $S_{p}$ be the $p$-local sphere. If $f : S_{p}\rightarrow B$ is any formal Zariski open
immersion then $L:=\mathbb{R}f_{*}\mathbb{L}f^{*}$ is clearly a smashing localization functor
in the sense of \cite[\S 3]{hps}. Its category $C$ of perfect\footnote{The word \textit{finite} instead of perfect would be more customary in this setting.} acyclics (i.e. perfect objects $X$ in $\mathrm{\mathbf{Ho}}(S_{p}-\mathbf{Mod})$ such that $LX$ is null)
is then a localizing thick subcategory of the homotopy category $\mathrm{\mathbf{Ho}}(S_{p}-\mathbf{Mod}^{\textrm{perf}})$ of the category of perfect $S_{p}$-modules, and therefore by \cite{hs} it is equivalent to the category $C_{n}$ of perfect $E(n)$-acyclics, for some $0\leq n < \infty$, where $E(n)$ is the $n$-th Johnson-Wilson $S_{p}$-module (see e.g. \cite{rav}); in other words $L$ and $L_{n}:=L_{E(n)}$ are both smashing localization functors on $\mathrm{\mathbf{Ho}}(S_{p}-\mathbf{Mod})$ having the same subcategory of finite acyclics. Therefore, if we
assume (one of the form of) the Telescope conjecture (see \cite{mil}), we get that $L_{n}$ and $L$
have equivalent categories of acyclics and so have equivalent categories of local objects.
But the category of local objects for $L$ is equivalent to the category
$\mathrm{\mathbf{Ho}}(B-\mathbf{Mod})$ (since $\mathbb{R}f_{*}$ is fully faithful by hypothesis) and
the category of local objects for $L_{n}$ is equivalent to the category $\mathrm{\mathbf{Ho}}((L_{n}S_{p})-\mathbf{Mod})$,
by \cite{wo} since $L_{n}$ is smashing. This easily implies that
the two commutative $S_{p}$-algebras $B$ and $L_{n}S_{p}$ are equivalent (i.e. isomorphic in $\mathbf{Ho}(S_{p}-\mathbf{Alg})$).\\
In conclusion, one sees that if
the Telescope conjecture is true, then, up to equivalence of $S_{p}$-algebras, the only (non-trivial) formal
Zariski open immersions for $S_{p}$ are given by the family
$$\mathcal{U}:=\left\{S_{p}\rightarrow L_{n}S_{p}\right\}_{0\leq n < \infty}.$$ This example shows that
the formal Zariski topology might be better suited in certain contexts than the
Zariski topology itself (e.g. it is not clear that there exists
any non-trivial Zariski open immersion of $S_{p}$, i.e. that the morphisms of commutative $S$-algebras
$S_{p} \longrightarrow L_{n}S_{p}$ are of finite presentation). Note however that the family $\mathcal{U}$ is not a formal Zariski \textit{covering} according to Definition \ref{d1} because the family of base-change functors  $$\left\{(-)\;\wedge^{\mathbb{L}}_{S_{p}} L_{n}S_{p}:\mathrm{\mathbf{Ho}}(S_{p}-Mod)\longrightarrow \mathrm{\mathbf{Ho}}(L_{n}S_{p}-Mod)\right\}_{0\leq n < \infty}$$ is not conservative; in fact, as Neil Strickland pointed out to us, the Brown-Comenetz dual $I$ of $S_{p}$ is a non-perfect non-trivial $S_{p}$-module which is nonetheless $L_{n}$-acyclic for any $n$. However, it is true that the family of base-changes above is conservative when restricted to the (homotopy) categories of perfect modules. Therefore, one could modify the second covering condition in Definition \ref{d1}, by only requiring the property of being conservative on the subcategories of perfect modules and relaxing the finiteness of $J$; let us call this modified covering condition \textit{formal Zariski covering-on-finites} condition. Then, $\mathcal{U}$ is a formal Zariski covering-on-finites family and indeed the unique one, up to equivalences of $S_{p}$-algebras, if the Telescope conjecture holds.
\item The previous example also shows that the commutative $S$-algebras
$L_{n}S_{p}$ are \textit{local for the formal Zariski topology} (again assuming
the Telescope conjecture). Indeed, for any
formal Zariski open covering $\{L_{n}S_{p} \longrightarrow B_{i}\}_{i\in I}$ there is
an $i$ such that $L_{n}S_{p} \longrightarrow B_{i}$ is an equivalence of commutative
$S$-algebras.
\end{enumerate}
}
\end{rmk}

\end{subsection}

\begin{subsection}{The brave new \'etale topology}

Notions of \'etale morphisms of commutative $S$-algebras has been studied
by several authors (\cite{ro, min2}). In this paragraph we present the definition
that appeared in \cite{hagI} and was used there in order to define the \'etale $K$-theory of
commutative $S$-algebras. \\

We refer to \cite{ba} for the notions of topological cotangent spectrum and of topological Andr\'e-Quillen cohomology relative to a morphism $A\rightarrow B$ of commutative $S$-algebras, except for slightly different notations. We denote by $\mathbb{L}\Omega_{B/A}\in\mathbf{Ho}(B-\mathbf{Mod})$ the topological cotangent spectrum (denoted as $\Omega_{B/A}$ in \cite{ba}) and, for any $B$-module $M,$ by $$\mathbb{L}Der_{A}(B,M):=\mathbb{R}\underline{Hom}_{A-\mathbf{Alg}/B}(B,B\vee M)$$ the derived space of topological derivations from $B$ to $M$ ($B\vee M$ being the trivial extension of $B$ by $M$). Note that there is an isomorphism $\mathbb{L}Der_{A}(B,M)\simeq \mathbb{R}\underline{Hom}_{B-\mathbf{Mod}}(\mathbb{L}\Omega_{B/A},M)$, natural in $M$.

\begin{df}\label{d2}
\begin{itemize}
\item
Let $f : A \longrightarrow B$ be a morphism of commutative $S$-algebras.
\begin{itemize}

\item The morphism $f$ is called \emph{formally \'etale} if
$\mathbb{L}\Omega_{B/A}\simeq 0$.

\item The morphism $f$ is called \emph{\'etale} if it
is formally \'etale and of finite presentation (as
a morphism of commutative $S$-algebras).
\end{itemize}

\item A family of morphisms $\{f_{i} : A \longrightarrow A_{i}\}_{i\in I}$
in $\mathbf{S-Alg}$ is called a
\emph{(formal) \'etale covering} if it satisfies the following two conditions.

\begin{itemize}

\item Each morphism $A \longrightarrow A_{i}$ is (formally) \'etale.

\item There exists a finite subset $J \subset I$ such that the
family of inverse image functors
$$\{\mathbb{L}f_{j}^{*} : \mathbf{Ho}(A-\mathbf{Mod}) \longrightarrow \mathbf{Ho}(\mathbf{A_{j}-Mod})\}_{j\in J}$$
is conservative (i.e. a morphism in $\mathbf{Ho}(A-\mathbf{Mod})$ is an isomorphism if and only if
its images by all the $\mathbb{L}f_{j}^{*}$'s are isomorphisms).
\end{itemize}
\end{itemize}
\end{df}

As shown in \cite[\S 5.2]{hagI}, (formal) \'etale covering families are stable by
equivalences, compositions and homotopy push-outs, and therefore
define a model topology on the model category $\mathbf{S-Aff}$.
Therefore one gets two model topologies called the \emph{brave new \'etale topology}
and the \emph{brave new formal \'etale topology}. The corresponding model sites
will be denoted by $(\mathbf{S-Aff},\textrm{\'et})$ and $(\mathbf{S-Aff},\textrm{f\'et})$, and will be called
the \emph{brave new \'etale site} and the \emph{brave new formal \'etale site}. \\

As for the brave new Zariski topology one proves that the brave new \'etale
topology is a generalization of the usual \'etale topology.

\begin{lem}\label{l4}
\begin{enumerate}
\item
Let $R \longrightarrow R'$ be a morphism of commutative rings.
The induced morphism $HR \longrightarrow HR'$ is an \'etale morphism
of commutative $S$-algebras (in the sense of Definition \ref{d2}) if and only if
the morphism $Spec\, R' \longrightarrow Spec\, R$ is an \'etale morphism of schemes.
\item
A family of morphisms of
commutative rings, $\{R \longrightarrow R'_{i}\}_{i\in I}$, induces
an \'etale covering family of commutative $S$-algebras $\{HR \longrightarrow HR_{i}'\}_{i\in I}$
(in the sense of Definition \ref{d2})
if and only if the family $\{Spec\, R_{i} \longrightarrow Spec\, R\}_{i\in I}$
is an \'etale covering of schemes.
\end{enumerate}
\end{lem}

\textit{Proof:} This is proved in \cite[\S 5.2]{hagI}. \hfill $\Box$ \\

Let $\mathbf{Aff}$ be the opposite category of commutative rings, and
$(\mathbf{Aff},\textrm{\'et})$ the big \'etale site. The site
$(\mathbf{Aff},\textrm{\'et})$ can also be considered as a model site (for the trivial
model structure on $\mathbf{Aff}$). Lemma \ref{l4} shows in particular that
the Eilenberg-Mac Lane functor $H : \mathbf{Aff} \longrightarrow \mathbf{S-Aff}$
induces a continuous morphism of model sites (\cite[\S 4.8]{hagI}). In this way, the
site $(\mathbf{Aff},\textrm{\'et})$ becomes a \emph{sub-model site}
of $(\mathbf{S-Aff},\textrm{\'et})$. \\

Another important fact is that the brave new \'etale topology is finer
than the brave new Zariski topology.

\begin{lem}\label{l5}
\begin{enumerate}
\item Any formal Zariski open immersion of commutative $S$-algebras
is a formally \'etale morphism.

\item Any Zariski open immersion of commutative $S$-algebras
is an \'etale morphism.

\item Any (formal) Zariski open covering of a commutative $S$-algebra
is a (formal) \'etale covering.

\end{enumerate}
\end{lem}

\textit{Proof:} Only $(1)$ requires a proof, and the proof will be similar
to the one of Lemma \ref{l1} $(2)$. Let $f : A\longrightarrow B$ be
a formal Zariski open immersion of commutative $S$-algebras. As the functor
$\mathbb{R}f_{*} : \mathbf{Ho}(B-\mathbf{Mod}) \longrightarrow \mathbf{Ho}(A-\mathbf{Mod})$
is a full embedding so is the induced functor
$\mathbb{R}f_{*} : \mathbf{Ho}(B-\mathbf{Alg}) \longrightarrow \mathbf{Ho}(A-\mathbf{Alg})$.
By definition of topological derivations one has for any $B$-module $M$,
$\mathbb{L}Der_{A}(B,M)=
\mathbb{R}\underline{Hom}_{A-\mathbf{Alg}/B}(B,B\vee M)$.
This and the fact that $\mathbb{R}f_{*}$ is fully faithful imply that
$$\mathbb{L}Der_{A}(B,M)=
\mathbb{R}\underline{Hom}_{A-\mathbf{Alg}/B}(B,B\vee M)\simeq
\mathbb{R}\underline{Hom}_{B-\mathbf{Alg}/B}(B,B\vee M)\simeq *,$$
and therefore that $\mathbb{L}\Omega_{B/A}\simeq 0$. \hfill $\Box$ \\

Lemma \ref{l5} implies that the identity functor of $\mathbf{S-Aff}$
defines a continuous morphism between model sites
$$(\mathbf{S-Aff},\mathrm{Zar}) \longrightarrow (\mathbf{S-Aff},\textrm{\'et}),$$
which is a base change functor from the brave new Zariski site to
the brave new \'etale site. The same is true for the formal versions of these sites. \\

To finish this part, we would like to mention a stronger version of
the brave new \'etale topology, called the \emph{$thh$-\'etale topology}, which is
sometimes more convenient to deal with.

\begin{df}\label{d3}
\begin{itemize}
\item
Let $f : A \longrightarrow B$ be a morphism of commutative $S$-algebras.
\begin{itemize}

\item The morphism $f$ is called \emph{formally $thh$-\'etale} if
for any commutative $A$-algebra $C$ the mapping
space $\mathbb{R}\underline{Hom}_{A-\mathbf{Alg}}(B,C)$
is $0$-truncated (i.e. equivalent to a discrete space).

\item The morphism $f$ is called \emph{$thh$-\'etale} if it
is formally $thh$-\'etale and of finite presentation (as
a morphism of commutative $S$-algebras).
\end{itemize}

\item A family of morphisms $\{f_{i} : A \longrightarrow A_{i}\}_{i\in I}$
in $\mathbf{S-Alg}$ is called a
\emph{(formal) $thh$-\'etale covering} if it satisfies the following two conditions.

\begin{itemize}

\item Each morphism $A \longrightarrow A_{i}$ is (formally) $thh$-\'etale.

\item There exists a finite subset $J \subset I$ such that the
family of inverse image functors
$$\{\mathbb{L}f_{j}^{*} : \mathbf{Ho}(A-\mathbf{Mod}) \longrightarrow \mathbf{Ho}(\mathbf{A_{j}-Mod})\}_{j\in J}$$
is conservative (i.e. a morphism in $\mathbf{Ho}(A-\mathbf{Mod})$ is an isomorphism if and only if
its images by all $\mathbb{L}f_{j}^{*}$ are isomorphisms).
\end{itemize}
\end{itemize}
\end{df}

It is easy to check that (formal) $thh$-\'etale coverings define a model topology on the
model category $\mathbf{S-Aff}$, call the \emph{(formal) $thh$-\'etale topology}.
The model category $\mathbf{S-Aff}$ together with these topologies will be
called the \emph{brave new $thh$-\'etale site} and the \emph{brave new formal $thh$-\'etale site}, denoted by
$(\mathbf{S-Aff},\textrm{thh-\'et})$ and $(\mathbf{S-Aff},\textrm{fthh-\'et})$, respectively. An equivalent way of stating the formal $thh$-\'etaleness condition for $A\rightarrow B$ is to say that the natural map $B\rightarrow S^{1}\otimes^{\mathbb{L}}B$ in $\mathbf{Ho}(A-\mathbf{Alg})$ is an isomorphism, or equivalently (by \cite{msv}), that the canonical map $B\rightarrow \mathrm{THH}(B/A,B)$ is an isomorphism in $\mathbf{Ho}(A-\mathbf{Alg})$, where $\mathrm{THH}$ denotes the topological Hochschild cohomology spectrum (see e.g. \cite[\S IX]{ekmm}). This equivalent characterization follows from the observation that, for any
morphism $f:B\rightarrow C$ of commutative $A$-algebras, one has an isomorphism (in $\mathbf{Ho}(\mathbf{SSet})$) of
mapping spaces:
$$\mathrm{Map}_{B-\mathbf{Alg}}(\mathrm{THH}(B/A,B),C)\simeq \Omega_{f}\mathrm{Map}_{A-\mathbf{Alg}}(B,C),$$
where $\Omega_{f}$ denotes the loop space at $f$. Therefore, the canonical map $B\rightarrow \mathrm{THH}(B/A,B)$ is an isomorphism in $\mathbf{Ho}(A-\mathbf{Alg})$ iff for any such $f$, $\Omega_{f}\mathrm{Map}_{A-\mathbf{Alg}}(B,C)$ is
contractible, i.e. iff the simplicial set $\mathrm{Map}_{A-\mathbf{Alg}}(B,C)\simeq \mathbb{R}\underline{Hom}_{A-\mathbf{Alg}}(B,C)$ is $0$-truncated. \\
This explains the name of this topology and, since as observed in \cite{min2} the Goodwillie derivative of $\mathrm{THH}$ is the suspension of the topological Andr\'e-Quillen spectrum $\mathrm{TAQ}$ (where, for any $B$-module $M$,  $\mathrm{TAQ}(B/A;M)$ is defined as the derived internal Hom from $\mathbb{L}\Omega_{B/A}$ to $M$ in the model category of $B$-modules) also shows that
(formal) $thh$-\'etale morphisms are (formal) \'etale morphisms. Therefore the identity functor
induces continuous morphisms of model sites
$$(\mathbf{S-Aff},\textrm{thh-\'et})\longrightarrow (\mathbf{S-Aff},\textrm{\'et}) \qquad
(\mathbf{S-Aff},\textrm{fthh-\'et})\longrightarrow (\mathbf{S-Aff},\textrm{f\'et}).$$
 We refer to
\cite{min2} for more details on the notion of $thh$-\'etale morphisms. \\

\end{subsection}

\begin{subsection}{Standard topologies}\label{stop}

Standard model topologies on $\mathbf{S-Aff}$ are obvious
extensions of usual Grothendieck topologies on affine schemes. They are
defined in the following way.

Let $\tau$ be one of the usual Grothendieck topologies on affine schemes
(i.e. Zariski, Nisnevich, \'etale or faithfully flat).

\begin{df}\label{dstandard}
A family of morphisms of commutative $S$-algebras
$\{A \longrightarrow B_{i}\}_{i\in I}$ is a \emph{standard $\tau$-covering} (also
called \emph{strong $\tau$-covering}) if it satisfies the following two
conditions.
\begin{itemize}
\item The induced family of morphisms of schemes $\{Spec\;\pi_{0}(B_{i}) \longrightarrow Spec \;\pi_{0}(A)\}_{i\in I}$
is a $\tau$-covering of affine schemes.

\item For any $i\in I$ the natural morphism of $\pi_{0}(B_{i})$-modules
$$\pi_{*}(A)\otimes_{\pi_{0}(A)}\pi_{0}(B_{i}) \longrightarrow
\pi_{*}(B_{i})$$
is an isomorphism.

\end{itemize}
\end{df}

Its easy to check that this defines a model topology $\tau^{s}$ on $\mathbf{S-Aff}$, called
the \emph{standard $\tau$-topology}.
The model site $(\mathbf{S-Aff},\tau^{s})$ may be called the \textit{brave new standard-$\tau$ site}.
The importance of standard topologies is that all $\tau^{s}$-coverings
of commutative $S$-algebras of the form $HR$ comes from usual
$\tau$-coverings of the scheme $Spec\, R$. Its behavior is therefore very close
to the geometric intuition one gets in Algebraic Geometry. \\

Finally, let us also mention the \emph{semi-standard} (or \emph{semi-strong})
model topologies. A family of morphisms of commutative $S$-algebras
$\{A \longrightarrow B_{i}\}_{i\in I}$ is a \emph{semi-standard $\tau$-covering} (also
called \emph{semi-strong $\tau$-covering}) if the induced
family of morphism of commutative graded rings $\{\pi_{*}(A) \longrightarrow \pi_{*}(B_{i})\}_{i\in I}$
is a $\tau$-covering. This also defines a model topology $\tau^{ss}$  on $\mathbf{S-Aff}$.

Both the standard and semi-standard type model sites (and $\mathbf{S}$-stacks over them,
see Section \ref{sstacks})  could be of some interest in the study of geometry over even,
periodic $S$-algebras (e.g. for elliptic spectra as in \cite{ahs}).

\end{subsection}

\end{section}

\begin{section}{$\mathbf{S}$-stacks and geometric $\mathbf{S}$-stacks}\label{sstacks}

Let $(M,\tau)$ be a \textit{model site} (i.e. a model category $M$ endowed with
a model topology $\tau$ in the sense of \cite{hagI}). Associated to it one
has a model category of \textit{prestacks} $M^{\wedge}$ and of
\textit{stacks} $M^{\sim,\tau}$. For details concerning these model categories we refer
to \cite[\S 4]{hagI}, and for the sake of brevity we only recall the following
facts.

\begin{itemize}

\item The model category $M^{\wedge}$ is a left Bousfield localization of
the model category $SSet^{M^{op}}$, of simplicial presheaves on
$M$ together with the projective levelwise model structure. The
local objects for this Bousfield localization are precisely the
simplicial presheaves $F : M^{op} \longrightarrow SSet$
which are equivalences preserving.

\item The model category $M^{\sim,\tau}$ is a left Bousfield localization (\cite[\S 3]{hi})
of the model category of prestacks $M^{\wedge}$, and the localization
(left Quillen) functor from $M^{\wedge}$ to $M^{\sim,\tau}$ preserves (up to equivalences) finite homotopy limits
(i.e. homotopy pull-backs).
The local objects for this Bousfield localization are the simplicial presheaves
$F : M^{op} \longrightarrow SSet$ which satisfy the following two conditions.
\begin{itemize}

\item The functor $F$ preserves equivalences (i.e. is a local object in $M^{\wedge}$).

\item For any $\tau$-hypercover $U_{*} \longrightarrow X$ in the model site
$(M,\tau)$ (\cite[\S 4.4]{hagI}), the induced morphism
$$F(X) \longrightarrow F(U_{*})$$
is an equivalence.

\end{itemize}
There is an \textit{associated stack} functor
$a:\mathbf{Ho}(M^{\wedge})\rightarrow \mathbf{Ho}(M^{\sim,\tau})$
right adjoint to the inclusion $\mathbf{Ho}(M^{\sim,\tau})\hookrightarrow \mathbf{Ho}(M^{\wedge})$.

\item There is a homotopical variant
\footnote{If $x\in M$, $\mathbb{R}\underline{h}(x)$ essentially sends $y\in M$ to
the mapping space $\mathrm{Map}_{M}(y,x).$} $$\mathbb{R}\underline{h}:
\mathbf{Ho}(M)\hookrightarrow \mathbf{Ho}(M^{\wedge})$$ of the \textit{Yoneda embedding} (\cite[\S 4.2]{hagI}).
\end{itemize}

Specializing to our present situation, where $M=\mathbf{S-Aff}$, we have one model category
$\mathbf{S-Aff}^{\wedge}$ of prestacks and zounds of model categories stacks
$$ \qquad
\mathbf{S-Aff}^{\sim,\mathrm{Zar}}, \qquad \mathbf{S-Aff}^{\sim,\textrm{\'et}}, \qquad
\mathbf{S-Aff}^{\sim,\textrm{thh-\'et}},$$
$$\mathbf{S-Aff}^{\sim,\textrm{fZar}}, \qquad \mathbf{S-Aff}^{\sim,\textrm{f\'et}}, \qquad
\mathbf{S-Aff}^{\sim,\textrm{fthh-\'et}},$$
$$\qquad \mathbf{S-Aff}^{\sim,\textrm{Zar}^{s}}, \qquad
\mathbf{S-Aff}^{\sim,\textrm{\'et}^{s}}, \qquad \mathbf{S-Aff}^{\sim,\textrm{ffqc}^{s}},$$
$$\dots \; etc \dots$$
These model categories come with right Quillen functors (the morphism of change of sites)
$$\xymatrix{
\mathbf{S-Aff}^{\sim, \textrm{\'et}} \ar[r] & \mathbf{S-Aff}^{\sim,\textrm{thh-\'et}} \ar[r] &
\mathbf{S-Aff}^{\sim,\textrm{Zar}} \ar[r] & \mathbf{S-Aff}^{\wedge}}$$
$$\xymatrix{
\mathbf{S-Aff}^{\sim,\textrm{f\'et}} \ar[r] & \mathbf{S-Aff}^{\sim,\textrm{fthh-\'et}} \ar[r] &
\mathbf{S-Aff}^{\sim,\textrm{fZar}} \ar[r] & \mathbf{S-Aff}^{\wedge}}$$
$$\xymatrix{
\mathbf{S-Aff}^{\sim,\textrm{\'et}} \ar[r] & \mathbf{S-Aff}^{\sim,\textrm{f\'et}} \ar[r] &
\mathbf{S-Aff}^{\sim,\textrm{fZar}} \ar[r] & \mathbf{S-Aff}^{\wedge}}$$
$$\xymatrix{
\mathbf{S-Aff}^{\sim,\textrm{\'et}^{s}} \ar[r] & \mathbf{S-Aff}^{\sim,\textrm{\'et}} \ar[r] &
\mathbf{S-Aff}^{\sim,\textrm{Zar}} \ar[r] & \mathbf{S-Aff}^{\wedge}}$$
$$\dots \; etc \dots$$
which allow to compare the various topologies on $\mathbf{S-Aff}$.

\begin{df}\label{d4}
Let $\tau$ be a model topology on the model site
$\mathbf{S-Aff}$.
\begin{itemize}
\item
The model category of \emph{$\mathbf{S}$-stacks for the topology $\tau$} is
$\mathbf{S-Aff}^{\sim,\tau}$.
\item A simplicial presheaf $F \in SPr(\mathbf{S-Aff})$ is called
an \emph{$S$-stack} if it is a local object in $\mathbf{S-Aff}^{\sim,\tau}$
(i.e. preserves equivalences and satisfies the descent property for
$\tau$-hypercovers).
\item
Objects in
the homotopy category $\mathbf{Ho}(\mathbf{S-Aff}^{\sim,\tau})$ will simply be called
\emph{$\mathbf{S}$-stacks} (without referring, unless it is necessary, to the underlying topology).
\end{itemize}
\end{df}

The category of $\mathbf{S}$-stacks, being the homotopy category of a model category,
has all kind of homotopy limits and colimits. Moreover, one can show that it has internal Hom's.
Actually, the model category of
$\mathbf{S}$-stacks is a \textit{model topos} in the sense of \cite[\S 3.8]{hagI}
(see also \cite{segtop}), and therefore behaves
very much in the same way as a category sheaves (but in a homotopical sense). In practice this
is very useful as it allows to use a lot of usual properties of simplicial
sets in the context of $\mathbf{S}$-stacks (in the same way as a lot of usual properties
of sets are true in any topos). \\

The Eilenberg-Mac Lane functor $H$ from commutative rings to commutative
$S$-algebras induces left Quillen functors
$$H_{!}: \mathbf{Aff}^{\sim,\tau_{0}} \longrightarrow \mathbf{S-Aff}^{\sim,\tau},$$
where $\tau_{0}$ is one of the standard topologies on affine schemes (e.g. $\textrm{Zar}$,
$\textrm{\'et}$, $\textrm{ffqc}, $\dots),
and $\tau$ is one of its possible extension to the model category $\mathbf{S-Aff}$
(e.g. $\textrm{Zar}$ can be extended to $\textrm{Zar}^{s}$ or to $\textrm{Zar}$, etc.). Here,
$\mathbf{Aff}^{\sim,\tau}$ is the usual model category of simplicial presheaves
on the Grothendieck site $(\mathbf{Aff},\tau)$ (with the projective
model structure \cite{bl}).
By deriving on the left one gets a functor
$$\mathbb{L}H_{!} : \mathbf{Ho}(\mathbf{Aff}^{\sim,\tau_{0}}) \longrightarrow \mathbf{Ho}(\mathbf{S-Aff}^{\sim,\tau}).$$
Therefore, our category of $\mathbf{S}$-stacks receives a functor from
the homotopy category of simplicial presheaves. In particular,
sheaves on affine schemes (and in particular
the category of schemes itself), and also $1$-truncated simplicial presheaves
(and in particular the homotopy category of algebraic stacks) can be all viewed as
examples of $\mathbf{S}$-stacks. However, one should be careful that the
functor $\mathbb{L}H_{!}$ has no reason to be fully faithful in general, though this is
the case for all the standard extensions (but not semi-standard) described in \S \ref{stop} (the reason for this is
that all covering families of some $HR$ are in fact induced from covering families of affine schemes. In
particular the restriction functor from $\mathbf{S-Aff}^{\sim,\tau} \longrightarrow \mathbf{Aff}^{\sim,\tau_{0}}$
will preserve local equivalences.). \\

If instead of requiring descent with respect to all $\tau$-hypercovers, we only require
descent with respect to those $\tau$-hypercovers which arise as homotopy nerves of $\tau$-covers,
we obtain the following weaker notion of stack.

\begin{df}\label{cechstacks}
Let $\tau$ be a model topology on $\mathbf{S-Aff}$. A simplicial presheaf $F:\mathbf{S-Aff}^{op}\rightarrow \mathbf{SSet}$ is said to be a \emph{$\check{\mathrm{C}}$ech} \emph{$\mathbf{S}$-stack} with respect to $\tau$ if it preserves equivalences and satisfies the following $\check{C}$ech descent condition. For any $\tau$-cover $\mathcal{U}=\left\{U_{i}\rightarrow X\right\}$, denoting by $\check{C}(\mathcal{U})_{*}$ its homotopy nerve,
the canonical map $$F(X)\longrightarrow \mathrm{holim}F(\check{C}(\mathcal{U})_{*})$$ is an isomorphism in $\mathbf{Ho}(\mathbf{SSet})$.
\end{df}

This weaker notion of stacks is the
one used recently by J. Lurie in \cite{lu} and has also appeared for stacks over Grothendieck sites in \cite{dhi}.\\

Note that, similarly to the case of $\mathbf{S}$-stacks, there is a model category $\mathbf{S-Aff}_{\check{C}}^{\sim,\tau}$ of $\check{\mathrm{C}}$ech $\mathbf{S}$-stacks, defined as the left
Bousfield localization of $\mathbf{S-Aff}^{\wedge}$ with respect to all the $\check{\mathrm{C}}$ech-nerves,
whose homotopy category is equivalent to the full subcategory of $\mathbf{Ho}(\mathbf{S-Aff}^{\wedge})$
consisting of $\check{\mathrm{C}}$ech $\mathbf{S}$-stacks.\\
In general the inclusion of $\mathbf{S}$-stacks into $\check{\mathrm{C}}$ech $\mathbf{S}$-stacks is proper;
however, one can prove (by adapting \cite[Prop. 6.1]{hs} to the context of model sites), that given any $n$-truncated
(equivalence preserving) simplicial presheaf $F$ on $\mathbf{S-Aff}$, $F$ is an $\mathbf{S}$-stack iff it is a $\check{\mathrm{C}}$ech $\mathbf{S}$-stack (regardless of the model topology $\tau$). The reader might wish to
read Appendix A of \cite{dhi} for more comparison results between stacks and $\check{\mathrm{C}}$ech stacks over usual Grothendieck sites.

\begin{subsection}{Some descent theory}\label{descent}

With the notations above, one can compose the Yoneda embedding
$\mathbb{R}\underline{h}:\mathbf{Ho}(\mathbf{S-Alg})^{op}\rightarrow \mathbf{Ho}(\mathbf{S-Aff}^{\wedge})$
with the associated stack functor $a:\mathbf{Ho}(\mathbf{S-Aff}^{\wedge})\rightarrow
\mathbf{Ho}(\mathbf{S-Aff}^{\sim,\tau})$ and obtain the derived $Spec$ functor
$$\mathbb{R}Spec\, : \mathbf{Ho}(\mathbf{S-Alg})^{op}=\mathbf{Ho}(\mathbf{S-Aff}) \longrightarrow
\mathbf{Ho}(\mathbf{S-Aff}^{\sim,\tau}),$$
for any model topology $\tau$ on $\mathbf{S-Aff}$.

\begin{df}\label{d5}
The topology $\tau$ is \emph{sub-canonical} (respectively, \emph{$\check{\mathrm{C}}$ech-subcanonical}) if for any $A\in \mathbf{S-Alg}$,
$\mathbb{R}\underline{h}_{A}$ is an $\mathbf{S}$-stack (resp., a $\check{\mathrm{C}}$ech $\mathbf{S}$-stack).
\end{df}

Note that $\tau$ is sub-canonical iff the functor
$\mathbb{R}Spec$ is fully faithful, and that sub-canonical implies $\check{\mathrm{C}}$ech-subcanonical. \\
Knowing whether a given model topology $\tau$ is sub-canonical or not is known as
the \emph{descent problem for $\tau$}, and in our opinion is a crucial question.
At present, we do not know if all the model topologies presented in the previous Section are
sub-canonical, and it might be that some of them are not.
The following lemma gives examples of sub-canonical topologies.

\begin{lem}\label{l6}
The (semi-)standard Zariski, Nisnevich, \'etale and flat model topologies of \S \ref{stop} are
all sub-canonical.
\end{lem}

\textit{Sketch of proof:} Let $\tau$ be one of these topologies,
$A$ be a commutative $S$-algebra, and $A \longrightarrow B_{*}$ be a
$\tau$-hypercover (\cite[\S 4.4]{hagI}). Using the fact that $\pi_{*}(B_{n})$
is flat over $\pi_{*}(A)$ for any $n$, one can check that
the cosimplicial $\pi_{*}(A)$-algebra $\pi_{*}(B_{*})$ is
again a $\tau$-hypercover of commutative rings. By usual
descent theory for affine schemes this implies that
the cohomology groups of the total complex of $[n] \mapsto \pi_{*}(B_{n})$
vanish except for $H^{0}(\pi_{*}(B_{*}))\simeq \pi_{*}(A)$. This implies that
the spectral sequence for the $\mathrm{holim}$
$$H^{p}([n] \mapsto \pi_{q}(B_{n})) \Rightarrow \pi_{p-q}(\mathrm{holim} B_{*})$$
degenerates at $E_{2}$ and that $\pi_{*}(A) \longrightarrow \pi_{*}(\mathrm{holim} B_{*})$
is an isomorphism. \hfill $\Box$ \\

Concerning the brave new Zariski topology one has the following partial result.

\begin{lem}\label{l7}
Let $\{A \longrightarrow A_{i}\}_{i\in I}$ be a finite Zariski covering family of
commutative $S$-algebras. Let $A \longrightarrow B=\vee_{i}A_{i}$ be the
coproduct morphism. Let $A \longrightarrow B_{*}$ be the cosimplicial
commutative $A$-algebra defined by
$$B_{n}:=\underbrace{B\wedge^{\mathbb{L}}_{A}B\wedge^{\mathbb{L}}_{A}\dots
\wedge^{\mathbb{L}}_{A}B}_{(n+1)\; times}$$
(i.e. homotopy co-nerve of the morphism $A \longrightarrow B$).
Then the induced morphism
$$A \longrightarrow \mathrm{holim}_{n\in \Delta}B_{n}$$
is an equivalence.
\end{lem}

\textit{Sketch of proof:} By definition of Zariski open immersion it
is not hard to see that the cosimplicial commutative $A$-algebra
$B_{*}$ is $m$-coskeletal, where $m$ is the cardinality of $I$. This means the
following: let $i_{m} : \Delta_{\leq m} \longrightarrow \Delta$ be the inclusion functor
form the full sub-category of objects $[i]$ with $i\leq m$. Then, one has
an equivalence of commutative $A$-algebras $B_{*}\simeq \mathbb{R}(i_{m})_{*}i_{m}^{*}(B_{*})$
(here $(i_{m}^{*},\mathbb{R}(i_{m})_{*})$ is the derived adjunction between
$\Delta$-diagrams and $\Delta_{\leq m}$-diagrams). From this one deduces easily that
$$\mathrm{holim}_{n\in \Delta}B_{n}\simeq \mathrm{holim}_{n\in \Delta_{\leq m}}B_{n}.$$
In particular, $\mathrm{holim}_{n\in \Delta}B_{n}$ is in fact a finite homotopy limit and therefore
will commute with the base change from $A$ to $B$, i.e.
$$(\mathrm{holim}_{n\in \Delta}B_{n})\wedge^{\mathbb{L}}_{A}B\simeq
\mathrm{holim}_{n\in\Delta}(B_{n}\wedge^{\mathbb{L}}B).$$
Now, as the functor
$\mathbf{Ho}(A-\mathbf{Mod})\longrightarrow \mathbf{Ho}(B-\mathbf{Mod})$ is
conservative (since the family $\{A \longrightarrow A_{i}\}_{i\in I}$
is a Zariski covering), one can replace $A$ by $B$ and the $A_{i}$
by $B\wedge_{A}^{\mathbb{L}}A_{i}$, and in particular one can
suppose that $A \longrightarrow B$ has a section. But, it is well known that
any morphism $A \longrightarrow B$ which has a section is such that
$A\simeq \mathrm{holim}_{n}B_{n}$ (the section can in fact be used in order to construct a retraction).
\hfill $\Box$ \\

\begin{cor}\label{zariscechsubcan}
The Zariski topology on $\mathbf{S-Aff}$ is $\check{C}$ech subcanonical.
\end{cor}

The results of the next sections, though stated for $\mathbf{S}$-stacks will also be correct by replacing
``$\mathbf{S}$-stack'' with ``$\check{\mathrm{C}}$ech $\mathbf{S}$-stack'', and
``subcanonical'' with ``$\check{\mathrm{C}}$ech subcanonical''.

\end{subsection}

\begin{subsection}{The $\mathbf{S}$-stack of perfect modules}\label{perfect}

Let $\tau$ be a model topology on $\mathbf{S-Aff}$.
One defines the $\mathbf{S}$-prestack $\underline{\mathsf{Perf}}$ of perfect modules in the following way.
For any commutative $S$-algebra $A$, we consider
the category $\mathsf{Perf}(A)$, whose objects are
perfect and cofibrant $A$-modules, and whose morphisms are equivalences
of $A$-modules. The pull back functors define a pseudo-functor
$$\begin{array}{cccc}
\mathsf{Perf} : & \mathbf{S-Alg} & \longrightarrow & \mathbf{Cat} \\
 & A & \longmapsto & \mathsf{Perf}(A) \\
 & (A \rightarrow B) & \longmapsto & (-\wedge_{A}B: \mathsf{Perf}(A)\rightarrow \mathsf{Perf}(B)).
\end{array}$$
Making this pseudo-functor into a strict functor from
$\mathbf{S-Alg}$ to $\mathbf{Cat}$ (\cite[Th. 3.4]{may}), and applying the
classifying space functor $\mathbf{Cat}\rightarrow \mathbf{SSet}$, we get
a simplicial presheaf denoted by $\underline{\mathsf{Perf}}$.

The following theorem relies on the so called strictification theorem (\cite[A.3.2]{hagI}), and its
proof will appear in \cite{hagII}.

\begin{thm}\label{t1}
The object $\underline{\mathsf{Perf}}$ is an $\mathbf{S}$-stack (i.e. satisfies the descent
condition for all $\tau$-hypercovers) iff the model topology $\tau$ is subcanonical.
\end{thm}

Another way to state Theorem \ref{t1} is by saying that $\tau$ is subcanonical iff, for any
commutative $S$-algebra $A$, the natural morphism
$$\underline{Hom}_{\mathbf{S-Aff}^{\sim,\tau}}(Spec\, A,\underline{\mathsf{Perf}})\simeq \underline{\mathsf{Perf}}(A) \longrightarrow
\mathbb{R}\underline{Hom}_{\mathbf{S-Aff}^{\sim,\tau}}(Spec\, A,\underline{\mathsf{Perf}})$$
is an equivalence of simplicial sets.

The $S$-stack of perfect complexes is a brave new analog of the stack of
vector bundles, and is of fundamental importance in brave new algebraic geometry.

\end{subsection}

\begin{subsection}{Geometric $\mathbf{S}$-stacks}\label{geom}

In this paragraph we will work with a fixed sub-canonical model topology $\tau$
on the model site $\mathbf{S-Aff}$.  We will define the notion of
\emph{geometric $\mathbf{S}$-stack}, which roughly speaking are \emph{quotients}
of affine $\mathbf{S}$-stacks by a \emph{smooth affine groupoid}. They will be brave new generalizations of Artin algebraic stacks (see \cite{lm}). In order to
state the precise definition, one first needs a notion
of smoothness for morphisms of commutative $S$-algebras. \\

For any perfect $S$-module $M$ one has
the (derived) free commutative $S$-algebra over $M$, $S \longrightarrow \mathbb{L}F_{S}(M)$.
For any commutative $S$-algebra $A$, one gets a morphism
$$A \longrightarrow A\wedge^{\mathbb{L}}_{S}\mathbb{L}F_{S}(M)\simeq
\mathbb{L}F_{A}(A\wedge^{\mathbb{L}}_{S}M).$$ Any morphism $A
\longrightarrow B$ in $\mathbf{Ho}(A-\mathbf{Alg})$ which is isomorphic to
 such a morphism will be called a \emph{perfect morphism} of
commutative $S$-algebras (and we will also say that $B$ is a
\emph{perfect commutative $A$-algebra}).

\begin{df}\label{d6}
A morphism of commutative $S$-algebras $f : A \longrightarrow B$ is
called \emph{smooth} if it satisfies the following two conditions.
\begin{itemize}
\item The $A$-algebra $B$ is finitely presented.

\item There exists an \'etale covering family $\{v_{i}: B\rightarrow B'_{i}\}_{i\in I}$ and,
for any $i\in I$, a homotopy commutative square
of commutative $S$-algebras
$$\xymatrix{A \ar[r]^-{f} \ar[d]_-{u} & B \ar[d]^-{v_{i}} \\
A' \ar[r]_-{f_{i}'} & B'_{i},}$$
where $f_{i}'$ is a perfect morphism,
and $u$ is an \'etale morphism.
\end{itemize}
\end{df}

One checks easily that smooth morphisms are stable by compositions
and homotopy base changes. Furthermore, any \'etale morphism is smooth, and therefore so
is any Zariski open immersion. \\

\textbf{Assumption:} At this point we will assume that the notion of smooth morphisms is \textit{local}
with respect to the chosen model topology $\tau$.  \\

This assumption will insure that
the notion of geometric $S$-stack, to be defined below,
behaves well.\\

Some terminology:

\begin{itemize}
\item
Let us come back to our homotopy category $\mathbf{Ho}(\mathbf{S-Aff^{\sim,\tau}})$ of $\mathbf{S}$-stacks, and
the Yoneda embedding (or derived $Spec$)
$$\mathbb{R}Spec : \mathbf{Ho}(\mathbf{S-Alg})^{op} \longrightarrow \mathbf{Ho}(\mathbf{S-Aff^{\sim,\tau}}).$$
The essential image of $\mathbb{R}Spec$ is called the category of \emph{affine $\mathbf{S}$-stacks}, which
is therefore anti-equivalent to the homotopy category of commutative $S$-algebras.
We will also call affine $\mathbf{S}$-stack any object in $\mathbf{S-Aff^{\sim,\tau}}$ whose
image in $\mathbf{Ho}(\mathbf{S-Aff^{\sim,\tau}})$ is an affine $\mathbf{S}$-stack. Clearly,
affine $\mathbf{S}$-stacks are stable by homotopy limits (indeed
$\mathrm{holim}_{i}(\mathbb{R}Spec\, A_{i})\simeq \mathbb{R}Spec\, (\mathrm{hocolim}_{i}A_{i})$).

\item A morphism of affine $\mathbf{S}$-stacks is called \textit{smooth} (over $S$) (resp. \'etale, a Zariski open immersions \dots)
if the corresponding one in $\mathbf{Ho}(\mathbf{S-Alg})$ is so.

\item A \emph{Segal groupoid object} in $\mathbf{S-Aff^{\sim,\tau}}$ is a simplicial object
$$X_{*} : \Delta^{op} \longrightarrow \mathbf{S-Aff^{\sim,\tau}}$$
which satisfies the following two conditions.
\begin{itemize}
\item For any $n\geq 1$, the $n$-th Segal morphism
$$X_{n} \longrightarrow \underbrace{X_{1}\times_{X_{0}}^{h}X_{1}\times^{h}_{X_{0}}\dots
X_{1}}_{n\; \mathrm{times}}$$
is an equivalence (in the model category $\mathbf{S-Aff^{\sim,\tau}}$ of $\mathbf{S}$-stacks).
When this condition is satisfied, it is well known that one can define a composition law (well defined
up to homotopy)
$$\mu : X_{1}\times_{X_{0}}^{h}X_{1} \longrightarrow X_{1}.$$
\item The induced morphism
$$(\mu,pr_{2}) : X_{1}\times_{X_{0}}^{h}X_{1} \longrightarrow X_{1}\times_{X_{0}}^{h}X_{1}$$
is an equivalence (i.e. the composition law is invertible up to homotopy).
\end{itemize}

\item For any simplicial object $X_{*} : \Delta^{op} \longrightarrow \mathbf{S-Aff^{\sim,\tau}}$,
we will denote by $|X_{*}|$ the homotopy colimit of $X_{*}$ in the model category
$\mathbf{S-Aff^{\sim,\tau}}$.

\end{itemize}

We are now ready to define geometric $\mathbf{S}$-stacks.

\begin{df}\label{d7}
An $\mathbf{S}$-stack $F$ is called \emph{geometric}
if it is equivalent to some $|X_{*}|$, where
$X_{*}$ is a Segal groupoid in $\mathbf{S-Aff^{\sim,\tau}}$ satisfying the
following two additional conditions.
\begin{itemize}

\item The $\mathbf{S}$-stacks $X_{0}$ and $X_{1}$ are affine $\mathbf{S}$-stacks.

\item The morphism $d_{0}: X_{1} \longrightarrow X_{0}$ is a smooth morphism
of affine $\mathbf{S}$-stacks.

\end{itemize}
\end{df}

The theory of geometric $\mathbf{S}$-stacks can then be pursued along the same lines as the theory of algebraic stacks (as done in \cite{lm}). For example, one can
define the notions of quasi-coherent and perfect modules on
a geometric $\mathbf{S}$-stack, $K$-theory of a geometric $\mathbf{S}$-stack (using perfect modules on it), higher geometric $\mathbf{S}$-stacks (such as
$2$-geometric $\mathbf{S}$-stacks), etc. We
refer the reader to \cite{hagII} for details. \\

We will finish this paragraph with the definition of the \textit{tangent} $\mathbf{S}$\textit{-stack}
and its main properties. \\

First of all, one defines a commutative $S$-algebra $S[\varepsilon]:=S\vee S$, which is the
trivial extension of $S$ by $S$. The $S$-algebra $S[\varepsilon]$ can be thought as
the \emph{brave new algebra of dual numbers}, i.e. the analog of $\mathbb{Z}[\varepsilon]$.
For any commutative $S$-algebra $A$, one has $A[\varepsilon]:=A\wedge^{\mathbb{L}}_{S}S[\varepsilon]\simeq A\vee A$, the
commutative $A$-algebra of dual numbers over $A$.

For any $\mathbf{S}$-stack $F \in \mathbf{S-Aff^{\sim,\tau}}$, one defines the \emph{tangent $\mathbf{S}$-stack of $F$}
as
$$\begin{array}{cccc}
TF: & \mathbf{S-Alg} & \longrightarrow & SSet \\
 & A & \mapsto  & F(A[\varepsilon]).
\end{array}$$

The tangent $\mathbf{S}$-stack $TF$ comes equipped with a natural projection
$p : TF \longrightarrow F$. One first notice that if $F$ is a geometric $\mathbf{S}$-stack (over any
base $A$), then
so is $TF$. Furthermore, the homotopy fibers of the projection $p$ are \emph{linear} $\mathbf{S}$-stacks in the following sense.
Let $A$ be a commutative $S$-algebra and $x : \mathbb{R}Spec\, A \longrightarrow F$
be a morphism of $\mathbf{S}$-stacks, i.e. an $A$-point of $F$. One considers
the homotopy pull back
$$\xymatrix{F_{x} \ar[r] \ar[d] & TF \ar[d] \\
\mathbb{R}Spec\, A \ar[r]_-{x} & F.}$$
Then, one can show that there exists an $A$-module $M$, such that
$F_{x}$ is equivalent (as a stack over $\mathbb{R}Spec\, A$) to
$\mathbb{R}Spec( \mathbb{L}F_{A}(M))$. In other words, one has a natural equivalence
$$F_{x}(B)\simeq \mathbb{R}\underline{Hom}_{A-\mathbf{Mod}}(M,B)$$
for any commutative $A$-algebra $B$. The $A$-module $M$ is called
the \emph{cotangent complex} of $F$ at the point $x$, and denoted by $\mathbb{L}\Omega_{F,x}$. Its derived dual $A$-module
$D(\mathbb{L}\Omega_{F,x})$ is called the \emph{tangent space} of $F$ at $x$.

\end{subsection}

\end{section}

\begin{section}{Some examples of geometric $\mathbf{S}$-stacks}

In this last Section we present two examples of
geometric $\mathbf{S}$-stacks. The first one arises from a classification
problem in Algebraic Topology, whereas the second
one is directly related to topological modular forms.

The first of
these examples (see \S \ref{mod}) shows that moduli
spaces in Algebraic Topology are not only discrete homotopy types
(as e.g. in \cite{bdg}), but have some additional rich
geometric structures very similar to the moduli spaces one studies
in Algebraic Geometry. It seems to us one of the most
simple non trivial example of brave new moduli stack,
and a relevant test for the whole theory.

The second example (see \S \ref{tmf})
seems to us much more intriguing, and might give new insights on
the construction and the study of topological modular forms.
We think that this research direction is definitely worth being investigated in the future, and
therefore we present a key open question that could be the starting point
of such an investigation. \\

We will work with a fixed subcanonical model topology $\tau$ on $\mathbf{S-Aff}$.

\begin{subsection}{The brave new group scheme $\mathbb{R}$\underline{Aut}(M) }

We fix a perfect $S$-module $M$, and we are going to define a
\emph{group $\mathbf{S}$-stack $\mathbb{R}\underline{Aut}(M)$}, of
auto-equivalences of $M$. This group $\mathbf{S}$-stack will be a
generalization of the group scheme $GL_{n}$, since it will be shown to
be an \emph{affine} and \emph{smooth} group $\mathbf{S}$-stack. Like many algebraic
stacks in Algebraic Geometry are quotients of affine schemes by $GL_{n}$, our example of
a geometric $\mathbf{S}$-stack in \S \ref{mod} will be a quotient of an affine $\mathbf{S}$-stack by
$\mathbb{R}\underline{Aut}(M)$ for
some $S$-module $M$. \\

For any commutative $S$-algebra $A$, one first defines
$$\mathbb{R}\underline{End}(M)(A):=
\mathbb{R}\underline{Hom}_{A-\mathbf{Mod}}(A\wedge^{\mathbb{L}}M,A\wedge^{\mathbb{L}}M).$$
Using for example the Dwyer-Kan simplicial localization techniques (\cite{dk1,dk2}), one can
make $A \mapsto \mathbb{R}\underline{End}(M)(A)$ into a functor
from $\mathbf{S-Alg}$ to the category $\mathbf{SMon}$ of simplicial monoids
$$\begin{array}{cccc}
\mathbb{R}\underline{End}(M) : & \mathbf{S-Alg} & \longrightarrow & \mathbf{SMon} \\
 & A & \longmapsto & \mathbb{R}\underline{End}(M)(A).
\end{array}$$

This defines $\mathbb{R}\underline{End}(M)$ as a monoid object
in $\mathbf{S-Aff^{\sim,\tau}}$. As its underlying object
in $\mathbf{S-Aff^{\sim,\tau}}$ is an $\mathbf{S}$-stack (for example using Theorem
\ref{t1}), we will say that $\mathbb{R}\underline{End}(M)$ is a \emph{monoid $\mathbf{S}$-stack}.

\begin{lem}\label{l10}
The $\mathbf{S}$-stack $\mathbb{R}\underline{End}(M)$ is affine and the structural morphism
$\mathbb{R}\underline{End}(M) \longrightarrow \mathbb{R}Spec\, S$
is perfect (hence smooth).
\end{lem}

\textit{Proof:} This is clear as
$$\mathbb{R}\underline{End}(M)\simeq
\mathbb{R}Spec\, (\mathbb{L}F_{S}(M\wedge_{S}^{\mathbb{L}}D(M))). $$ \hfill $\Box$ \\

For any commutative $S$-algebra $A$, one defines
$\mathbb{R}\underline{Aut}(M)(A)$ to be the sub-monoid
of $\mathbb{R}\underline{End}(M)(A)$ consisting of auto-equivalences. In other words, $\mathbb{R}\underline{Aut}(M)(A)$ is defined by the following homotopy pull-back diagram in $\mathbf{SSet}$ $$\xymatrix{ \mathbb{R}\underline{Aut}(M)(A)\ar[r] \ar[d] & \mathbb{R}\underline{End}(M)(A) \ar[d]\\
[M\wedge^{\mathbb{L}}A,M\wedge^{\mathbb{L}}A]' \ar@{^{(}->}[r] & [M\wedge^{\mathbb{L}}A,M\wedge^{\mathbb{L}}A]}$$ where $[-,-]'$ is the subset of isomorphisms in $\mathbf{Ho}(\mathbf{SSet})$.  This defines
a functor
$$\begin{array}{cccc}
\mathbb{R}\underline{Aut}(M) : & \mathbf{S-Alg} & \longrightarrow & \mathbf{SMon} \\
 & A & \longmapsto & \mathbb{R}\underline{Aut}(M)(A).
\end{array}$$
Once again the underlying object in $\mathbf{S-Aff^{\sim,\tau}}$
is an $\mathbf{S}$-stack, and therefore $\mathbb{R}\underline{Aut}(M)$ is a monoid
$\mathbf{S}$-stack. Furthermore, the monoid law on $\mathbb{R}\underline{Aut}(M)$
is invertible up to homotopy, and we will therefore say
that $\mathbb{R}\underline{Aut}(M)$ is a \emph{group $\mathbf{S}$-stack}.

\begin{lem}\label{l11}
The $\mathbf{S}$-stack $\mathbb{R}\underline{Aut}(M)$ is affine and the structural morphism
$\mathbb{R}\underline{Aut}(M) \longrightarrow \mathbb{R}Spec\, S$
is smooth. In other words, $\mathbb{R}\underline{Aut}(M)$ is
an \emph{affine} and \emph{smooth} group $\mathbf{S}$-stack.
\end{lem}

\textit{Proof:}
The following proof is inspired by the proof of \cite[I.9.6.4]{egaI}. Let $B$ be the commutative $S$-algebra $\mathbb{L}F_{S}(M\wedge_{S}^{\mathbb{L}}D(M))$ corresponding
to the affine $\mathbf{S}$-stack $\mathbb{R}\underline{End}(M)$. There exists a universal
endomorphism of $B$-modules
$$u : M\wedge^{\mathbb{L}}_{S}B \longrightarrow M\wedge^{\mathbb{L}}_{S}B$$
such that for any commutative $B$-algebra $C$,
the endomorphism
$$u\wedge_{B}^{\mathbb{L}} id_{C}: M\wedge^{\mathbb{L}}_{S}C \longrightarrow M\wedge^{\mathbb{L}}_{S}C$$
is equal (in $\mathbf{Ho}(B-\mathbf{Mod})$) to the corresponding point in
$$\mathbb{R}\underline{End}(M)(C)\simeq \mathbb{R}\underline{Hom}_{\mathbf{S-Alg}}(B,C).$$

Consider now the homotopy cofiber $K \in \mathbf{Ho}(B-\mathbf{Mod})$
of the universal endomorphism $u$. Clearly, $K$ is a perfect $B$-module, and one can
therefore consider the open Zariski immersion $B \longrightarrow B_{K}$ (Lemma \ref{l2}).
It is easy to check by construction that
$$\mathbb{R}\underline{Aut}(M) \simeq \mathbb{R}Spec\, B_{K},$$
which proves that $\mathbb{R}\underline{Aut}(M)$ is an affine $\mathbf{S}$-stack.
Finally, as the morphism $\mathbb{R}Spec\, B_{K} \longrightarrow
\mathbb{R}Spec\, B$ is smooth (being a Zariski open immersion), one sees
(using the fact that $\mathbb{R}Spec\, B$ is perfect hence smooth) that
$\mathbb{R}\underline{Aut}(M) \longrightarrow \mathbb{R}Spec\, S$ is also smooth. \hfill $\Box$ \\

\end{subsection}

\begin{subsection}{Moduli of algebra structures}\label{mod}

In this paragraph, we fix a perfect $S$-module $M$.
We will define an $\mathbf{S}$-stack $\underline{Ass}_{M}$, classifying
associative and unital algebras whose underlying module is $M$.

For any commutative $S$-algebra $A$, we have the category
$A-\mathbf{Ass}$, of associative and unital $A$-algebras (i.e. associative monoids in the monoidal category $(A-\mathbf{Mod},\wedge_{A})$);
these are new versions of the old $A_{\infty}$-ring spectra.
The category $A-\mathbf{Ass}$ has a model category structure
for which fibrations and equivalences are detected on the underlying objects
in $A-\mathbf{Mod}$. We denote by
$A-\mathbf{Ass}_{M}^{\mathrm{cof}}$ the subcategory of
$A-\mathbf{Ass}$ whose objects are cofibrant objects $B$ such that
there exists a $\tau$-covering family $\{A \longrightarrow A_{i}\}_{i\in I}$
such that each $A_{i}$-module $B\wedge^{\mathbb{L}}_{A}A_{i}$ is
equivalent to $M\wedge_{A}^{\mathbb{L}}A_{i}$ (we say that
the underlying $A$-module of $B$ is \emph{$\tau$-locally
equivalent to $M$}), and whose morphisms are equivalences of $A$-algebras.
The base change functors define a lax functor
$$\begin{array}{cccc}
Ass_{M}: & \mathbf{S-Alg} & \longrightarrow & \mathbf{Cat} \\
 & A & \mapsto & A-\mathbf{Ass}_{M}^{\mathrm{cof}} \\
 & (A\rightarrow B) & \mapsto & -\wedge_{A}B.
\end{array}$$
Strictifying this functor (\cite[Th. 3.4]{may}) and then applying the classifying space functor, one gets
a simplicial presheaf
$$\begin{array}{cccc}
\underline{Ass}_{M}: & \mathbf{S-Alg} & \longrightarrow & \mathbf{SSet} \\
 & A & \mapsto & B(A-\mathbf{Ass}_{M}^{\mathrm{cof}}).
\end{array}$$

For the following theorem, let us recall that for any commutative $S$-algebra $A$, any
associative and unital $A$-algebra $B$ and any $B$-bimodule $M$, one has
an $A$-module of $A$-derivations $Der_{A}(B,M)$ from $B$ to $M$.
This can be derived on the left (in the model category of
associative and unital $A$-algebras !) to
$\mathbb{L}Der_{A}(B,M)$.

\begin{thm}\label{t2}
Let $\tau$ be a subcanonical model topology on $\mathbf{S-Aff}$, and $M$ be a perfect $S$-module.
\begin{enumerate}
\item The object $\underline{Ass}_{M} \in \mathbf{S-Aff^{\sim,\tau}}$ is
an $\mathbf{S}$-stack.

\item The $\mathbf{S}$-stack $\underline{Ass}_{M}$ is geometric.


\end{enumerate}
\end{thm}

\textit{Sketch of proof:} Point $(1)$ can be proved with the same
techniques used in Theorem \ref{t1} and will not be proved here. We refer
to \cite{hagII} for details.


Let us prove part $(2)$ which is in fact a corollary of one of the main result
of C. Rezk thesis \cite{re}.

Let us start by considering the full sub-$\mathbf{S}$-stack of
$\underline{\mathsf{Perf}}$ (see \S \ref{perfect}) consisting of perfect modules which are $\tau$-locally equivalent to
$M$. By the result of Dwyer and Kan \cite[2.3]{dk3}, this $\mathbf{S}$-stack is clearly equivalent as an object
in $\mathbf{S-Aff^{\sim,\tau}}$ to $B\mathbb{R}\underline{Aut}(M)$, the classifying
simplicial presheaf of the group $\mathbf{S}$-stack $\mathbb{R}\underline{Aut}(M)$.
Forgetting the algebra structure gives
a morphism of $\mathbf{S}$-stacks
$$f : \underline{Ass}_{M} \longrightarrow B\mathbb{R}\underline{Aut}(M).$$
Using the techniques of equivariant stacks developed in
\cite{kpt} (or more precisely their straightforward extensions to the present context
of $\mathbf{S}$-stacks), one sees that the $\mathbf{S}$-stack
$\underline{Ass}_{M}$ is equivalent to the quotient $\mathbf{S}$-stack
$$[X/\mathbb{R}\underline{Aut}(M)],$$ where $X$ is the homotopy fiber of the morphism $f$ and
$\mathbb{R}\underline{Aut}(M)$ acts on $X$. By Lemma \ref{l11},  $\mathbb{R}\underline{Aut}(M)$ is an affine smooth group $\mathbf{S}$-stack, so we only need to show that $X$ is an affine $\mathbf{S}$-stack
(because the classifying Segal groupoid for the action of $\mathbb{R}\underline{Aut}(M)$ on $X$ will then
satisfies the conditions of Definition \ref{d7}).

Using \cite[Thm. 1.1.5]{re}, one sees that the homotopy fiber $X$
is equivalent to the $\mathbf{S}$-stack $$\mathbb{R}\underline{Hom}_{\mathrm{Oper}}(\mathcal{ASS},\underline{\mathcal{E}nd}(M)): A\longmapsto \mathbb{R}\underline{Hom}_{\mathrm{Oper}}(\mathcal{ASS},\mathcal{E}nd(M\wedge_{S}^{\mathbb{L}}A)),$$
where $\mathbb{R}\underline{Hom}_{\mathrm{Oper}}(\mathcal{ASS},\mathcal{E}nd(M\wedge_{S}^{\mathbb{L}}A))$ is the  derived Hom (or mapping space) of unital operad morphisms from the final operad $\mathcal{ASS}$ (classifying
associative and unital algebras) to the endomorphisms operad $\mathcal{E}nd(M\wedge_{S}^{\mathbb{L}}A)$
of the $A$-module $M\wedge_{S}^{\mathbb{L}}A$
(here operads are in the symmetric monoidal category $\mathbf{S}$ of $S$-modules).
This means that, for any commutative $S$-algebra $A$,
there is an equivalence
$$X(A)\simeq \mathbb{R}\underline{Hom}_{\mathrm{Oper}}(\mathcal{ASS},\underline{\mathcal{E}nd}(M))(A),$$
functorial in $A$.

Now, writing the operad $\mathcal{ASS}$ as a homotopy colimit
$$\mathcal{ASS}\simeq \mathrm{hocolim}_{n\in \Delta^{op}}\mathcal{O}_{n},$$
where each $\mathcal{O}_{n}$ is a free operad, one sees that
$$X\simeq \mathrm{holim}_{n\in \Delta}\mathbb{R}\underline{Hom}_{\mathrm{Oper}}(\mathcal{O}_{n},\underline{\mathcal{E}nd}(M)).$$
Since affine $\mathbf{S}$-stacks are stable under homotopy limits, it is therefore enough to check that the $\mathbf{S}$-stack $\mathbb{R}\underline{Hom}_{\mathrm{Oper}}(\mathcal{O},\underline{\mathcal{E}nd}(M))$
is affine for any free operad $\mathcal{O}$. But, saying that an operad
$\mathcal{O}$ is free means that there is a family $\left\{P_{m}\right\}_{m>0}$ of $S$-modules, and functorial (in $A\in S-\mathbf{Alg}$) equivalences
$$\mathbb{R}\underline{Hom}_{\mathrm{Oper}}(\mathcal{O},\underline{\mathcal{E}nd}(M))(A)\simeq
\prod_{m}\mathbb{R}\underline{Hom}_{S-\mathbf{Mod}}
(P_{m}\wedge_{S}^{\mathbb{L}}M^{\wedge^{\mathbb{L}} m}\wedge_{S}^{\mathbb{L}}D(M),A),$$ where the funny notation
$M^{\wedge^{\mathbb{L}} m}$ stands for the derived smash product $M\wedge^{\mathbb{L}}\dots\wedge^{\mathbb{L}}M$ of $M$ with itself $m$ times. So it is enough to show that for any $S$-module $P$, the (pre)stack $$A\longmapsto \mathbb{R}\underline{Hom}_{S-\mathbf{Mod}}
(P\wedge_{S}^{\mathbb{L}}M^{\wedge^{\mathbb{L}} m}\wedge_{S}^{\mathbb{L}}D(M),A)$$ is affine.
But this is clear since this stack is equivalent to $\mathbb{R}Spec\, B$ where $B$ is the derived free commutative $S$-algebra
$$B:=\mathbb{L}F_{S}(P\wedge_{S}^{\mathbb{L}}M^{\wedge^{\mathbb{L}} n}\wedge_{S}^{\mathbb{L}}D(M)).$$
This implies that $X$ is an affine $\mathbf{S}$-stack and completes the proof.
\hfill $\Box$ \\

\begin{rmk} \emph{The relationship between the tangent space of
$\underline{Ass}_{M}$ at a point $x : \mathbb{R}Spec\, A \longrightarrow \underline{Ass}_{M}$ (corresponding to
an associative $A$-algebra whose underlying $A$-module
is $\tau$-locally equivalent to $M\wedge_{S}^{\mathbb{L}}A$) and
the suspension $\mathbb{L}Der_{A}(B,B)[1]$ of the $A$-module of
derived $A$-derivations of the associative $A$-algebra $B$ into the $B$-bimodule $B$,
will be investigated elsewhere.}
\end{rmk}

Theorem \ref{t2} has also generalizations when one consider
algebra structures over a given operad (for example commutative algebra structures).
It can also be enhanced by considering
categorical structures such as $A_{\infty}$-categorical structures, as explained in \cite{hagdag};
the corresponding moduli space
gives an example of a $2$\textit{-geometric} $\mathbf{S}$-stack.\\

\end{subsection}

\begin{subsection}{Topological modular forms and geometric $\textbf{S}$-stacks}\label{tmf}

In this final section we will be working with
the standard \'etale topology on $\mathbf{S-Aff}$. The corresponding model
site will be denoted by $(\mathbf{S-Aff},\textrm{\'et}^{s})$, and its model category of
stacks by $\mathbf{S-Aff}^{\sim,\textrm{\'et}^{s}}$. We recall that the
topology $\textrm{\'et}^{s}$ is known to be subcanonical (see Lem. \ref{l6}).

As explained right after Def. \ref{d4}, the Eilenberg-MacLane spectrum construction
gives rise to a fully faithful functor
$$\mathbb{L}H_{!} : \mathbf{Ho}(\mathbf{Aff}^{\sim,\textrm{\'et}}) \longrightarrow \mathbf{Ho}(\mathbf{S-Aff}^{\sim,\textrm{\'et}^{s}}),$$
where the left hand side is the homotopy category of simplicial presheaves on the
usual \'etale site of affine schemes. This functor
has a right adjoint, called the \emph{truncation functor}
$$h^{0}:=H^{*} : \mathbf{Ho}(\mathbf{S-Aff}^{\sim,\textrm{\'et}^{s}}) \longrightarrow \mathbf{Ho}(\mathbf{Aff}^{\sim,\textrm{\'et}}),$$
simply given by composing a simplicial presheaf $F : \mathbf{S-Aff}^{op} \longrightarrow \mathbf{SSet}$
with the functor $H : \mathbf{Aff} \longrightarrow \mathbf{S-Aff}$. \\

Let us denote by $\overline{\mathcal{E}}$ the moduli stack of
generalized elliptic curves with integral geometric fibers,
which is the standard
compactification of the moduli stack of elliptic curves by adding
the nodal curves at infinity (see e.g. \cite[IV]{dera}, where it is denoted by $\mathcal{M}_{(1)}$); recall that
$\overline{\mathcal{E}}$ is a Deligne-Mumford stack, proper and smooth over $Spec\; \mathbb{Z}$ (\cite[Prop. 2.2]{dera}).

As shown by
recent works of M. Hopkins, H. Miller, P. Goerss,
N. Strickland, C. Rezk and M. Ando, there exists a natural
presheaf of commutative $S$-algebras on the small
\'etale site $\overline{\mathcal{E}}_{\textrm{\'et}}$ of
$\overline{\mathcal{E}}$. We will denote this presheaf by
$\mbox{\em tmf}$. Recall that by construction, if $U=Spec\, A \longrightarrow \overline{\mathcal{E}}$
is an \'etale morphism, corresponding to an elliptic curve $E$ over
the ring $A$, then $\mbox{\em tmf}(U)$ is the (connective) elliptic cohomology theory
associated to the formal group of $E$ (in particular, one
has $\pi_{0}(\mbox{\em tmf}(U))=A$). Recall also that the (derived) global sections
$\mathbb{R}\Gamma(\overline{\mathcal{E}},\mbox{\em tmf})$, form a
commutative $S$-algebra, well defined in $\mathbf{Ho}(\mathbf{S-Alg})$,
called the \emph{spectrum of topological modular forms}, and denoted by $\mathsf{tmf}$.

Let $U \longrightarrow \overline{\mathcal{E}}$ be a surjective \'etale morphism with
$U$ an affine scheme, and let us consider its nerve
$$\begin{array}{cccc}
U_{*} : & \Delta^{op} & \longrightarrow & \mathbf{Aff} \\
 & [n] & \mapsto & U_{n}:=\underbrace{U\times_{\overline{\mathcal{E}}}U\times_{\overline{\mathcal{E}}}\dots
\times_{\overline{\mathcal{E}}}U}_{n\; times}.
\end{array}$$

This is a simplicial object in $\overline{\mathcal{E}}_{\textrm{\'et}}$, and by
applying $\mbox{\em tmf}$ we obtain a
co-simplicial object in $\mathbf{S-Alg}$
$$\begin{array}{cccc}
\mbox{\em tmf}(U_{*}) : & \Delta & \longrightarrow & \mathbf{S-Alg} \\
 & [n] & \mapsto & \mbox{\em tmf}(U_{n}).
\end{array}$$
Taking $\mathbb{R}Spec$ (\S \ref{descent}) of this diagram we obtain
a simplicial object in the model category $\mathbf{S-Aff}^{\sim,\textrm{\'et}^{s}}$
$$\begin{array}{cccc}
\mathbb{R}Spec\; (\mbox{\em tmf}(U_{*})) : & \Delta^{op} & \longrightarrow & \mathbf{S-Aff}^{\sim,\textrm{\'et}^{s}} \\
 & [n] & \mapsto & \mathbb{R}Spec\; (\mbox{\em tmf}(U_{n})).
\end{array}$$

The homotopy colimit of this diagram will be denoted by
$$\overline{\mathcal{E}}_{\mathbf{S}}:=\mathrm{hocolim}_{n\in \Delta^{op}}
\mathbb{R}Spec\; (\mbox{\em tmf}(U_{*})) \in \mathbf{Ho}(\mathbf{S-Aff}^{\sim,\textrm{\'et}^{s}}).$$

The following result is just a remark as there is
essentially nothing to prove; however, we prefer to state it
as a theorem to emphasize its importance.

\begin{thm}\label{t3}
The stack $\overline{\mathcal{E}}_{\mathbf{S}}$ defined above is
a geometric $\mathbf{S}$-stack. Furthermore, there exists a natural
isomorphism in $\mathbf{Ho}(\mathbf{Aff}^{\sim,\textrm{\'et}})$
$$h^{0}(\overline{\mathcal{E}}_{\mathbf{S}})\simeq \overline{\mathcal{E}}.$$
\end{thm}

\textit{Proof:} To prove that $\overline{\mathcal{E}}_{\mathbf{S}}$ is geometric, it is enough to check that
the simplicial object $\mathbb{R}Spec\; (\mbox{\em tmf}(U_{*}))$ is a Segal groupoid
satisfying the conditions of Def. \ref{d7}. For this, recall that for
any morphism $U=Spec\, B \rightarrow V=Spec\, A$ in $\overline{\mathcal{E}}_{\textrm{\'et}}$, the natural morphism
$$\pi_{*}(\mbox{\em tmf}(V))\otimes_{\pi_{0}(\mbox{\em tmf}(V))}\pi_{0}(\mbox{\em tmf}(U))\simeq \pi_{*}(\mbox{\em tmf}(V))\otimes_{A}B \longrightarrow
\pi_{*}(\mbox{\em tmf}(U))$$
is an isomorphism. This shows that the functor
$$\mathbb{R}Spec\; (\mbox{\em tmf}(-)) : \overline{\mathcal{E}}_{\textrm{\'et}} \longrightarrow
\mathbf{S-Aff}^{\sim,\textrm{\'et}^{s}}$$
preserves homotopy fiber products and therefore sends Segal groupoid objects
to Segal groupoid objects. This shows in particular that
$\mathbb{R}Spec\; (\mbox{\em tmf}(U_{*}))$ is a Segal groupoid object. The same fact also shows that
for any morphism $U=Spec\, B \rightarrow V=Spec\, A$ in $\overline{\mathcal{E}}_{\textrm{\'et}}$,
the induced map $\mbox{\em tmf}(V) \longrightarrow \mbox{\em tmf}(U)$ is
a strong \'etale morphism in the sense of Def. \ref{dstandard}, and
therefore is an \'etale and thus smooth morphism. This implies that
$\mathbb{R}Spec\; (\mbox{\em tmf}(U_{*}))$ satisfies the conditions of \ref{d7} and therefore
shows that $\overline{\mathcal{E}}_{\mathbf{S}}$ is indeed a geometric
$\mathbf{S}$-stack.

The truncation functor $h^{0}$ clearly commutes with homotopy colimits, and therefore
$$h^{0}(\overline{\mathcal{E}}_{\mathbf{S}})\simeq
\mathrm{hocolim}_{n\in \Delta^{op}}h^{0}(\mathbb{R}Spec\; (\mbox{\em tmf}(U_{n}))) \in \mathbf{Ho}(\mathbf{Aff}^{\sim,\textrm{\'et}}).$$
Furthermore, for any connective affine $\mathbf{S}$-stack $\mathbb{R}Spec\, A$ one
has a natural isomorphism $h^{0}(\mathbb{R}Spec\, A)\simeq Spec\, \pi_{0}(A)$.
Therefore, one sees immediately that there is a natural isomorphism of
simplicial objects in $\mathbf{Aff}^{\sim,\textrm{\'et}}$
$$h^{0}(\mathbb{R}Spec\; (\mbox{\em tmf}(U_{*})))\simeq U_{*}.$$
Therefore, we get
$$h^{0}(\overline{\mathcal{E}}_{\mathbf{S}})\simeq
\mathrm{hocolim}_{n\in \Delta^{op}}h^{0}(\mathbb{R}Spec\; (\mbox{\em tmf}(U_{n})))\simeq
\mathrm{hocolim}_{n\in \Delta^{op}}U_{n}\simeq \overline{\mathcal{E}},$$
as $U_{*}$ is the nerve of an \'etale covering of $\overline{\mathcal{E}}$. \hfill $\Box$ \\

\begin{rmk}
\emph{The proof of Theorem \ref{t3} shows not only  
that $\overline{\mathcal{E}}_{\mathbf{S}}$ is geometric, but also
that it is a} strong Deligne-Mumford geometric \textbf{S}-stack, \emph{in the
sense that one can replace the word} smooth \emph{by} strongly \'etale \emph{in 
Definition \ref{d7}. See \cite{hagII} for further details.} 
\end{rmk}

Theorem \ref{t3} tells us that the presheaf of topological modular forms
$\mbox{\em tmf}$ provides a natural geometric $\mathbf{S}$-stack $\overline{\mathcal{E}}_{\mathbf{S}}$
whose truncation is the usual stack of elliptic curves
$\overline{\mathcal{E}}$. Furthermore, as the small strong \'etale topoi
of $\overline{\mathcal{E}}_{\mathbf{S}}$ and $\overline{\mathcal{E}}$
coincides (this is a general fact about strong \'etale model topologies, see \cite{hagII}), we see that
$$\mathsf{tmf}:=\mathbb{R}\Gamma(\overline{\mathcal{E}},\mbox{\em tmf}) \simeq
\mathbb{R}\Gamma(\overline{\mathcal{E}}_{\mathbf{S}},\mathcal{O}),$$
and therefore that topological modular forms can be simply interpreted
as \emph{functions on the geometric $\mathbf{S}$-stack $\overline{\mathcal{E}}_{\mathbf{S}}$}.
Of course, our construction of $\overline{\mathcal{E}}_{\mathbf{S}}$ has essentially been done to
make this true, so this is not a surprise. However, we have gained a bit from the conceptual
poit of view: since after all $\overline{\mathcal{E}}$ is a moduli stack, now
that we know the existence of the geometric $\mathbf{S}$-stack
$\overline{\mathcal{E}}_{\mathbf{S}}$ we can ask for a \textit{modular interpretation}
of it, or in other words for a direct geometric description of the corresponding
simplicial presheaf on $\mathbf{S-Aff}$. An answer to this question not only would
provide a direct construction
of $\mathsf{tmf}$, but would also give a conceptual interpretation of it in a geometric
language closer the usual notion of modular forms.

\begin{Q}\label{q1}
Find a modular interpretation of the
$\mathbf{S}$-stack
$\overline{\mathcal{E}}_{\mathbf{S}}$.
\end{Q}

Essentially, we ask which are the brave new ``objects'' that the
$\mathbf{S}$-stack $\overline{\mathcal{E}}_{\mathbf{S}}$ classifies.
Of course we do not know the answer to this question, though some
progress are being made by J. Lurie and the authors.
It seems that the \textbf{S}-stack $\overline{\mathcal{E}}_{\mathbf{S}}$ itself is not really
the right object to look at, and one should rather consider the non-connective
version of it (defined using the non-connective version of $\mbox{\em tmf}$) for
which a modular interpretation seems much more accessible. Though this modular interpretation
is still conjectural and not completely achieved, it does
use some very interesting notions of \emph{brave new abelian varieties},
\emph{brave new formal groups} and
their geometry. We think that the achievement of such a program could be
the starting point of a rather new and deep interaction between
stable homotopy theory and algebraic geometry, involving
many new questions and objects, but probably also
new insights on classical objects of algebraic topology.

\end{subsection}

\end{section}

\end{document}